\documentclass[12pt]{article}

\usepackage{amsfonts}
\usepackage{amssymb}
\usepackage{latexsym}
\usepackage{amsmath}
\usepackage{bibmods}

\usepackage{graphicx}
\usepackage{color}
\usepackage{bibnames}

\usepackage{amssymb}
\usepackage{amscd} %% Package for commutative diagrams
\usepackage{color}

 \newtheorem{thm}{Theorem}[section]
 \newtheorem{prop}{Proposition}[section]
 \newtheorem{lem}{Lemma}[section]
 \newtheorem{cor}{Corollary}[section]
\newtheorem{defn}{}[section]

 \newtheorem{rem}{Remark}[section]
\newcommand{\rbx}{\hfill{\rule{1ex}{1ex}}}

\newcommand{\tp}{\overset{\!\!\!\!\!\circ}{\sT_\alpha^+}}
\newcommand{\ind}{\mathrm{ind}\,}
\newcommand{\Span}{\mathrm{span}\,}
\newcommand{\diag}{\mathrm{diag}\,}
\newcommand{\wind}{\mathrm{wind}\,}
 
 \newcommand{\vp}{\varphi}
 \newcommand{\ve}{\varepsilon}
\newcommand{\D}{\displaystyle}

\newcommand{\Ind}{\mbox{\rm Ind }}
\newcommand{\smb}{\mbox{\rm smb}\,}
\newcommand{\coker}{\mbox{\rm coker}\,}

\newcommand{\im}{\mathrm{im}\,}

\newcommand{\nn}{\nonumber}

\newcommand{\re}{\textrm{Re}\,}

 \newcommand{\esssup}{\hbox{\rm ess}\sup_{\!\!\!\!\!\!\!\!\!\!\! t\in {\mathbb{T}}}}

\newcommand{\cB}{\mathcal{B}}

\newcommand{\cL}{\mathcal{L}}

\newcommand{\sC}{{\mathbb C}}

\newcommand{\sN}{{\mathbb N}}

\newcommand{\sR}{{\mathbb R}}
\newcommand{\sS}{{\mathbb S}}
\newcommand{\sT}{{\mathbb T}}

\newcommand{\sZ}{{\mathbb Z}}

\begin{document}

%%%%%%%%%%%%%%%%%%%%%%%%%%%%%%%%%%
%%%%%%%%%%%%%%%%%%%%%%%%%%%%%%
\vspace*{10mm}

\begin{center}
{\Large\textbf{Generalized Toeplitz plus Hankel operators: kernel
structure and defect numbers}\footnote{This work was partially
supported by the German Academic Exchange Service (DAAD) and by the
Universiti Brunei Darussalam, Grants UBD/GSR/S\&T/19 and
UBD/PNC2/2/RG/1(159)}}
\end{center}

\vspace{5mm}

%\author{Victor D. Didenko}
\begin{center}

\textbf{Victor D. Didenko and Bernd Silbermann}

 %\title[\textit{Invertibility of Toeplitz plus Hankel operators}]{}

\vspace{2mm}


%    author one information
Universiti Brunei Darussalam, Bandar Seri Begawan, BE1410  Brunei;
diviol@gmail.com

Technische Universit{\"a}t Chemnitz, Fakult{\"a}t f\"ur Mathematik,
09107 Chemnitz, Germany; silbermn@mathematik.tu-chemnitz.de

 \end{center}

%\keywords{Toeplitz plus Hankel operator, Invertibility}

  \vspace{10mm}

\textbf{2010 Mathematics Subject Classification:} Primary 47B35,
47B38; Secondary 47B33,  45E10

\textbf{Key Words:} Generalized Toeplitz plus Hankel operator,
Defect numbers, Kernel

%%%%%%%%%%%%%%%%%%%%%%%%%%%%%%%%%%%%
%%%%%%%%%%%%%%%%%%%%%%%%%%%%%%%%%%%

\begin{abstract}
Generalized Toeplitz plus Hankel operators $T(a)+H_{\alpha}(b)$
generated by functions $a,b$ and a linear fractional Carleman shift
$\alpha$ changing the orientation of the unit circle $\sT$ are
considered on the Hardy spaces $H^p(\sT)$, $1<p<\infty$. If the
functions $a,b\in L^\infty(\sT)$ and satisfy the condition
 $$
a(t) a(\alpha(t))=b(t) b(\alpha(t)),\quad t\in \sT,
 $$
the defect numbers of the operators $T(a)+H_{\alpha}(b)$ are
established and an explicit description of the structure of the
kernels and cokernels of the operators mentioned is given.
\end{abstract}

\section{Introduction\label{s1}}
Let $\sT:=\{ t\in\sC:|t|=1 \}$ be the counterclockwise oriented unit
circle in the complex plane $\sC$. By $PC:=PC(\sT)$ we denote the
set of all piecewise continuous function on $\sT$, i.e. the set of
all functions $f$ such that for any point $t_0\in \sT$ there are
finite left and right limits $f(t_0-0)$ and $f(t_0+0)$. As usual,
$C:=C(\sT)$ refers to the set of all continuous functions on $\sT$.
Further, let $L^p=L^p(\sT)$, $1\leq p\leq \infty$ stand for the
space of all Lebesgue measurable functions $f$ such that
 \begin{align*}
||f||_p &: = \left ( \int_\sT |f(t)|^p \,| dt| \right )^{1/p}
<\infty, 1\leq p <\infty;\\
 ||f||_\infty&: = \esssup |f(t)|<\infty.
\end{align*}
If $f\in L^1$, then by $\widehat{f}_n$ we denote the Fourier
coefficients of the function $f$,
 \begin{equation*}
\widehat{f}_n:=\frac{1}{2\pi}\int_0^{2\pi}f(e^{i\theta})
e^{-in\theta}\, d\theta\, ,
 \end{equation*}
and for $1\leq p \leq \infty$ let  $H^p:=H^p(\sT)$ and
$\overline{H^p}:=\overline{H^p(\sT)}$ be the Hardy spaces,
\begin{equation}\label{eq1}
 \begin{aligned}
   H^p:=\{ f\in L^p: \widehat{f}_n=0 \,\text{ for all } n <0\}, \\
 \overline{H^p}:=\{ f\in L^p: \widehat{f}_n=0 \,\text{ for all } n
 >0\}.
\end{aligned}
\end{equation}
On the space $L^p$, $1<p<\infty$, consider the operators $J$, $P$,
and $Q$ defined by
 \begin{equation}\label{eq3}
\begin{aligned}
 J:& f(t) \to t^{-1}f(t^{-1})\, ,\\
 P: & \sum_{n\in \sZ} \widehat{f}_n t^n \to \sum_{n\in \sZ_+}^\infty
 \widehat{f}_n t^n \, ,\\
  & Q:=I-P \, ,
  \end{aligned}
\end{equation}
where $I$ is the identity operator, and $\sZ_+:=\sN\cup \{0\}$ is
the set of all non-negative integers. The operators $P$ and $Q$ are
complimentary projections and it is well-known that they are bounded
on any space $L^p$ for $p\in(1,\infty)$, but not on the spaces $L^1$
and $L^\infty$.

 The classical Toeplitz plus Hankel operators $T(a)+H(b):H^p\to
 H^p$, $a,b\in L^\infty$,
  \begin{equation}\label{eq2}
    T(a)+H(b):= PaP + PbQJ
\end{equation}
are subject to numerous studies, where their Fredholm properties and
other related problems have been investigated. Recall that an
operator $A$ from the Banach algebra of all linear continuous
operators $\cB(X)$ acting on a Banach space $X$ is called Fredholm
if the range  $\im A:=\{y\in X: y=Ax, x\in X \}$ of $A$ is a closed
subspace of $X$ and the dimensions $\dim \ker A$ and $\dim \coker A$
of the subspaces
 \begin{align*}%\label{eq}
\ker A:=\{ x:X: Ax=0\}, \quad
 \coker A: =\{x\in X: A^*x=0 \}
\end{align*}
are finite. Here, $A^*$ denotes the adjoint operator to the operator
$A$. As far as the operator $T(a)+H(b)$ is concerned, for $a,b\in
PC$ its Fredholm properties can be immediately derived by a direct
applications of results \cite[Sections 4.95-4.102]{BS:2006},
\cite[Sections 4.5 and 5.7]{RSS:2011}, \cite{RS1990}. The case of
quasi piecewise continuous generating functions has been studied in
\cite{Si:1987}, whereas formulas for the index of the operators
\eqref{eq2} considered on various Banach and Hilbert spaces and with
various assumptions about the generating functions $a$ and $b$ have
been established in \cite{DS:2013,RS:2012}. Recently, progress has
been made in computation of defect numbers $\dim\ker (T(a)+H(b))$
and $\dim\coker (T(a)+H(b))$ for various classes of generating
functions $a$ and $b$ \cite{BE:2013,DS:2014a}. Moreover, a more
delicate problem of the description of the spaces $\ker (T(a)+H(b))$
and $\coker (T(a)+H(b))$ has been considered
\cite{DS:2014a,DS:2014}.

The aim of the present work is to study generalized Toeplitz plus
Hankel operators. These operators are similar to the classical
Toeplitz plus Hankel operator \eqref{eq3} but the flip operator $J$
of \eqref{eq2} is replaced by another operator $J_\alpha$ generated
by a linear fractional shift $\alpha$ changing the orientation of
the circle $\sT$. Areas of particular interest to us are the kernels
and cokernels of such operators and we are going to derive an
explicit description of these spaces in the case where the
generating functions $a$ and $b$ belong to the space $L^\infty$ and
satisfy an additional algebraic relation. Note that generalized
Toeplitz plus Hankel operators have been previously considered in
\cite{KLR:2007, KLR:2009} but under more restrictive assumptions.
Moreover, in the present work we single out certain classes of
generalized Toeplitz plus Hankel operators which are subject to
Coburn--Simonenko theorem. Recall that the classical
Coburn--Simonenko Theorem claims that if $a$ is a non-zero function,
then at least one of the numbers $\dim \ker T(a)$ or $\dim \coker
T(a)$ is equal to zero. For a Fredholm Toeplitz operator $T(a)$, the
generating function $a$ is invertible. Therefore, the
Coburn--Simonenko Theorem indicates that any Fredholm Toeplitz
operator is one-sided invertible.

Let $\beta$ be a complex number such that $|\beta|>1$, and let
$\sS^2$ denote the Riemann sphere. Consider the mapping
$\alpha:\sS^2 \to \sS^2$ defined by
 \begin{equation}\label{eq4}
   \alpha (z):=\frac{z-\beta}{\overline{\beta}z-1},
\end{equation}
and recall some basic properties of $\alpha$. In particular, one
has:

 \begin{enumerate}
\item The mapping $\alpha:\sS^2\to \sS^2$ is an one-to-one,
$\alpha(\sT)=\sT$,  and if $D^+:=\{ z\in\sC:|z|<1\}$ is the interior
of the unite circle $\sT$ and $\overline{D^+}:=D^+\cup \sT$, then
\begin{equation}\label{eq5}
\begin{aligned}
\alpha(D^+)=\sS^2\setminus \overline{D^+}, \\
\alpha(\sS^2\setminus \overline{D^+})=D^+ .
\end{aligned}
\end{equation}
It is clear that $\alpha$ is an automorphism of the Riemann sphere.
and the mappings $H^p \to \overline{H^p}$, $h\mapsto h\circ \alpha$
and $\overline{H^p} \to H^p$, $h\mapsto h\circ \alpha$ are
well-defined isomorphisms. Note that in the last relations, the
mapping $\alpha$ is understood as acting on the unit circle $\sT$. A
proof of this result can be given by using relations \eqref{eq9}. We
omit the details here but mention that they can be found in the
proof of Proposition 2.2 of \cite{KLR:2009}.

\item  The mapping $\alpha:\sT\to \sT$ changes the orientation
of $\sT$, satisfies the Carleman condition $\alpha(\alpha(t))=t$ for
all $t\in\sT$, and  possesses two fixed points $t_+=(1+
\lambda)/\overline{\beta}$ and $t_-=(1- \lambda)/\overline{\beta}$,
where $\lambda:=i \sqrt{|\beta|^2-1}$.

\item The mapping $\alpha$ admits the factorization
 \begin{equation}\label{eq6}
\alpha(t)=\alpha_+(t) \, t^{-1} \, \alpha_-(t),
\end{equation}
where
\begin{equation}\label{eq7}
\alpha_+(t)=\frac{t-\beta}{\lambda}, \quad \alpha_-(t)=
\frac{\lambda t}{\overline{\beta} t-1},
\end{equation}
and $\alpha_+^{\pm 1}\in H^\infty$, $\alpha_-^{\pm 1}\in
\overline{H^\infty}$.

\item On the space $L^p$, $1<p<\infty$, the mapping $\alpha$ generates
a bounded linear operator $J_\alpha$, called weighted shift operator
and defined by
  \begin{equation}\label{eq8}
J_\alpha \vp(t):= t^{-1} \alpha_-(t) \vp(\alpha(t)), \quad t\in
 \sT.
 \end{equation}
\end{enumerate}

Further, for $a\in L^\infty$ let $a_\alpha$ denote the composition
of the functions $a$ and $\alpha$, i.e.
 $$
 a_\alpha(t):=a(\alpha(t)), \quad t\in\sT.
 $$
It is clear that for any $a\in L^\infty$, the operator $J_{\alpha}$
satisfies the relation
  \begin{equation}\label{eq9.1}
  J_\alpha a I=a_\alpha J_\alpha,
 \end{equation}
and for any $n\in \sZ$, one has $(a^n)_{\alpha}=(a_{\alpha})^n:
=a_{\alpha}^n$.

In addition,
 \begin{gather}\label{eq9}
  J_\alpha^2 =I, \quad J_\alpha P=Q J_\alpha, \quad J_\alpha Q=
 PJ_\alpha,\\
\intertext{and}
 \alpha_+^{\pm1}(\alpha(t))=\alpha_-^{\pm1}(t), \quad
\overline{\alpha}_+^{\pm1}(t)=\alpha_-^{\mp1}(t). \label{eq9.3}
 \end{gather}
Indeed, consider for example the second identity in \eqref{eq9.3}.
Then
 $$
\overline{\alpha_+(t)} =\overline{\left
(\frac{t-\beta}{\lambda}\right )}=\frac{1/t-\overline{\beta}}{-
\lambda}
 =\frac{\overline{\beta}t-1}{t \lambda}=\alpha_-^{-1}(t).
 $$

Let us briefly discuss the situation with adjoint operators. Recall
that the adjoint to $H^p$ is the space $L^q/ \im Q$,
$p^{-1}+q^{-1}=1$. It can be identified with the space $H^q$ but
equipped with an equivalent norm coinciding with usual one in the
case $p=2$ (see \cite[Section 1.42]{BS:2006}). For our purposes
here, there is no need to distinguish between the two spaces $H^q$
and $L^q/ \im Q$, and we can also assume that the adjoint operator
$A^*$ to the operator $A\in \cL(H^P)$ is acting on the space $H^q$.

On the space $H^p$, $1<p<\infty$, any element $a\in L^\infty$
defines two operators  $T(a)$ and $H_\alpha(a)$ such that
 \begin{align*}%\label{eq}
   T(a)&: \vp \mapsto Pa \vp,\\
H_\alpha(a)&: \vp\mapsto PbQJ_\alpha \vp.
\end{align*}
The operator $T(a)$ is referred to as Toeplitz operator generated by
the function $a$, whereas the operator $H_\alpha(a)$ is called
generalized Hankel operator generated by the function $a$ and the
shift $\alpha$.  Many properties of generalized Hankel operators are
similar to the corresponding properties of classical Hankel
operators $PbQJ$. For example, Toeplitz and generalized Hankel
operators are connected in the following way,
  \begin{equation}\label{eq10}
\begin{aligned}
   T(cd)=T(c)T(d)+H_\alpha(c)H_\alpha (d_\alpha), \\
    H_\alpha(cd)=T(c)H_\alpha(d)+H_\alpha(c)T(d_\alpha).
\end{aligned}
\end{equation}
For classical Hankel operators such kind formulas are well known
\cite{BS:2006}, and they are often used in various studies of
Toeplitz and Hankel operators.

Let $a,b\in L^\infty$. In the present paper we study the properties
of the operators $T(a)+H_{\alpha}(b):H^p \to H^p$, $1<p<\infty$.
Such operators are called generalized Toeplitz plus Hankel operators
generated by the functions $a$, $b$ and by the shift $\alpha$ or
simply generalized Toeplitz plus Hankel operators. Similarly,
$T(a)-H_\alpha(b)$ are called generalized Toeplitz minus Hankel
operators.

Recall that the kernel and cokernel dimensions of generalized
Toeplitz plus Hankel operators and even more general operators have
been studied in \cite{KLR:2007, KLR:2009}. The approach of
\cite{KLR:2007, KLR:2009} is based on a special factorization of the
operators in question, which in turn involves factorizations of
matrix functions. Thus assuming that the corresponding operators are
Fredholm, the authors express the kernel and cokernel dimensions in
terms of the partial indices of some matrix valued functions.
However, it is well known that for an arbitrary matrix function the
computation of its partial indices is a very demanding task.
Moreover, in the most cases this is not possible at all. Therefore,
an important problem is to find certain classes of operators where
more efficient results can be derived. One such a class of
generalized Toeplitz plus Hankel operators $T(a)+H_{\alpha}(b)$,
$a,b\in L^\infty$ is mentioned in \cite{KLR:2009}. It can be
characterized by the relation
 $$
a a_{\alpha}= b b_{\alpha},
 $$
and in \cite{KLR:2009} such operators are studied under the
following assumptions.
\begin{enumerate}
    \item The operators $T(a)+H_{\alpha}(b):H^p\to H^p$ and
    $T(a)-H_{\alpha}(b):H^p\to H^p$ are Fredholm.
    \item The functions $a$ and $ba^{-1}$ admit \emph{bounded} Wiener--Hopf
    factorizations.
\end{enumerate}

In the present work, the conditions imposed on the operator
$T(a)+H_{\alpha}(b)$ are more general and a completely different
approach to these operators is used. In particular, condition i) is
replaced by a weaker one and we also avoid \emph{bounded}
Wiener--Hopf factorization. Recall that a Toeplitz operator $T(a)$
acting on $H^p$\/ is Fredholm if and only if the generating function
$a$ admits a Wiener--Hopf factorization in the sense of the
definition given in Section~\ref{s4}. The reader can also find there
a brief discussion of the two types of Wiener--Hopf factorizations
mentioned.

 This paper is organized as follows. In Section \ref{s2},
connections between the kernels of the operators $T(a)\pm
H_{\alpha}(b)$ and the kernel of a matrix Toeplitz operator
$T(V(a,b))$ are considered. Here we also describe the structure of
$\ker T(V(a,b))$ in terms of the kernels of certain scalar Toeplitz
operators. This allows us to establish a general representation for
the kernels of the operators $T(a)\pm H_{\alpha}(b)$. In Section
\ref{s3}, some classes of generalized Toeplitz plus Hankel
operators, where Coburn--Simonenko theorem holds, are studied. For
classical Toeplitz plus Hankel operators similar problems are
discussed in \cite{BE:2013, DS:2014a, DS:2014}, whereas
\cite{DS:2014b} deals with Wiener--Hopf plus Hankel operators. In
Section \ref{s4}, a decomposition for the kernels of special scalar
Toeplitz operators is derived. This decomposition is employed in
Section \ref{s5} in order to give a complete and effective
description for the kernels and cokernels of the generalized
Toeplitz plus Hankel operators under consideration. The remaining
part of the paper is devoted to the operators with piecewise
continuous generating functions $a$ and $b$.

\section{Kernels of generalized Toeplitz plus Hankel operators and related matrix
Toeplitz operators. \label{s2}}

There are well known relations between classical Toeplitz plus
Hankel operators and matrix Toeplitz operators. Similar formulas can
be also derived for the operators $T(a)+H_{\alpha}(b)$. Moreover,
these formulas can be effectively used in order to establish
connections between the kernels of the corresponding operators and
to derive a complete and efficient description of the kernels of the
operators under consideration. More precisely,  let $a,b\in
L^\infty$ and $a\in GL^\infty$, where $GL^\infty$ denotes the group
of invertible elements in $L^\infty$. The last assumption about the
generating function $a$ is justified by the fact that the
semi-Fredholmness of the operator $T(a)+H_{\alpha}(b)$ implies that
$a\in GL^\infty$. Along with $T(a)+H_{\alpha}(b)$ consider the
generalized Toeplitz minus Hankel operators $T(a)-H_{\alpha}(b)$ and
matrix Toeplitz operator $T(V(a,b))$, where
\begin{equation*}
V(a,b):=
\left(%
\begin{array}{cc}
 a- bb_{\alpha}a_{\alpha}^{-1} & d \\
 -c   &  a_{\alpha}^{-1}\\
   \end{array}%
\right), \quad c:=b_{\alpha}a_{\alpha}^{-1},\;d:=ba_{\alpha}^{-1}.
\end{equation*}

 The following result plays an important role in the description of
 the kernels of Toeplitz plus Hankel operators.

    \begin{lem}\label{l1}
Assume that $a,b\in L^\infty$, $a\in GL^\infty$, and the operators
$T(a)\pm H_{\alpha}(b)$ are considered on the space $H^p$,
$1<p<\infty$. Then
  \begin{itemize}
    \item If $(\vp,\psi)^T\in \ker T(V(a,b))$, then
  \begin{align}\label{eq11}
     (\Phi, \Psi)^T= &
      \frac{1}{2}(\vp-J_{\alpha}Qc\vp+J_{\alpha}Qa_{\alpha}^{-1}\psi,
      \vp+J_{\alpha}Qc\vp-J_{\alpha}Qa_{\alpha}^{-1}\psi)^T
      \\
      &\in
     \ker \diag(T(a)+H_{\alpha}(b),T(a)-H_{\alpha}(b)) &\nn
\end{align}

    \item If $(\Phi, \Psi)^T\in \ker\diag (T(a)+H_{\alpha}(b),T(a)-H_{\alpha}(b)
 )$, then
\begin{equation}\label{eq12}
(\Phi+\Psi, P(b_{\alpha}(\Phi + \Psi)
+a_{\alpha}J_{\alpha}P(\Phi-\Psi))^T\in \ker T(V(a,b)).
\end{equation}
    \end{itemize}
        Moreover, the operators
 \begin{align*}%\label{eq}
    \mathbf{U}_1: \ker T(V(a,b)) \to \ker\diag(T(a)+H_{\alpha}(b),T(a)-H_{\alpha}(b)),\\[1ex]
    \mathbf{U}_2: \ker\diag(T(a)+H_{\alpha}(b),T(a)-H_{\alpha}(b)) \to \ker T(V(a,b)),
\end{align*}
defined, respectively, by the relations \eqref{eq11} and
\eqref{eq12} are mutually inverses to each other.
 \end{lem}
\textbf{Proof.} First all all we note that the operator $T(a)+
H_{\alpha}(b)\in\cL(H^p(\sT))$ and $(T(a)+
H_{\alpha}(b))P+Q\in\cL(L^p(\sT))$ are simultaneously invertible or
Fredholm and have the same index. For simplicity, in what follows we
will write $T(a)\pm H_{\alpha}(b)+Q$ instead of $(T(a)\pm
H_{\alpha}(b))P+Q$. On the space $L^p\times L^p$ consider the
operator $\diag (T(a)+H_{\alpha}(b)+Q, T(a)-H_{\alpha}(b)+Q)$ and
represent it as the product of three matrix operators, viz.
\begin{equation}\label{eq13}
 \left(%
\begin{array}{cc}
  T(a)+H_{\alpha}(b)+Q &  0\\
   0 &  T(a)-H_{\alpha}(b)+Q\\
   \end{array}%
\right) =
 B(T(V(a,b))\!+\!\diag (Q,Q))A\,,
\end{equation}
where $A,B: L^p\times L^p\to L^p\times L^p$ are invertible
operators,
\begin{equation}\label{eq14}
A= \left(%
\begin{array}{cc}
  I &  0\\
  b_{\alpha}I  &  a_{\alpha}I\\
   \end{array}%
\right)
\left(%
\begin{array}{cc}
  I & I\\
  J_{\alpha}  &  -J_{\alpha}\\
   \end{array}%
\right),
\end{equation}
and $B=B_1 B_2 B_3$ with
 \begin{align*}
B_1 &=2\left(%
\begin{array}{cc}
  I & J_{\alpha}\\
  I  &  -J_{\alpha}\\
   \end{array}%
\right), \\
B_2&= \diag(I,I) -\diag(P,Q) \left(%
\begin{array}{cc}
  aI &  bI\\
  b_{\alpha}I  &  a_{\alpha}I\\
   \end{array}%
\right) \diag(Q,P), \\
B_3 &=\diag(I,I) +\diag(P,P) \left(%
\begin{array}{cc}
  a-b b_{\alpha} a_{\alpha}^{-1} &  d I\\
  -cI  &  a_{\alpha}^{-1}I\\
   \end{array}%
\right) \diag(Q,Q).
 \end{align*}
The representation \eqref{eq13} can be verified by straightforward
computations using relations \eqref{eq9.1}-\eqref{eq9}. The rest of
the proof goes through as for Lemma~3.2 of \cite{DS:2014}.
  \rbx
\sloppy
   \begin{rem}\label{r2}
Relation \eqref{eq13} shows that the matrix operators
$\diag(T(a)+H_{\alpha}(b),T(a)-H_{\alpha}(b))$ and $T(V(a,b))$ are
simultaneously Fredholm. This is not always true for the operators
$T(a)+H_{\alpha}(b)$ and $T(a)-H_{\alpha}(b)$ themselves. Moreover,
even if both operators $T(a)+H_{\alpha}(b)$ and $T(a)-H_{\alpha}(b)$
are Fredholm, they can have different indices. Nevertheless, the
relation \eqref{eq13} plays an important role in what follows and
can be employed in order to obtain a description of the kernels of
generalized Toeplitz plus Hankel operators.
 \end{rem}

From now on we will always assume that $a$ belongs to the group
$GL^\infty$ of invertible elements from $L^\infty$ and that the
generating functions $a$ and $b$ satisfy the condition
\begin{equation}\label{eq15}
   a a_{\alpha}= b b_{\alpha}.
\end{equation}
Relation \eqref{eq15} is called matching condition, whereas the
corresponding duo $(a,b)$ is referred to as an $\alpha$-matching
pair or simply a matching pair. For each $\alpha$-matching pair
$(a,b)$ one can assign another $\alpha$-matching pair $(c,d)$, where
$c=ab^{-1}$ and $d= b a_{\alpha}^{-1}$. Such a pair $(c,d)$ is
called the subordinated pair for $(a,b)$, and it is easily seen that
the functions which constitutes a subordinated pair have a specific
property, namely $cc_{\alpha}=1=dd_{\alpha}$. Throughout this paper
any function $g\in L^\infty$ satisfying the condition
 $$
g g_\alpha=1
 $$
is called $\alpha$-matching or simply matching function. In passing
note that for the subordinated pair $(c,d)$, functions $c$ and $d$
can also be expressed in the form
\begin{equation*}
c=b_{\alpha}a_{\alpha}^{-1}, \quad d= b_{\alpha}^{-1} a.
\end{equation*}
Besides, if $(c,d)$ is the  subordinated  pair for an
$\alpha$-matching pair $(a,b)$, then $(\overline{d},\overline{c})$
is the subordinated pair for the matching pair $(\overline{a},
\overline{b}_\alpha)$ defining the adjoint operator
  \begin{equation}\label{adjoint}
(T(a)+H_{\alpha}(b))^*=T(\overline{a})+H_{\alpha}(\overline{b}_\alpha),
 \end{equation}
to the operator $T(a)+H_{\alpha}(b)$. Note that in order to derive
formula \eqref{adjoint}, one can use the relation
   \begin{equation}\label{eq9.2}
   J_{\alpha}^* \psi (t)= t^{-1} \alpha_-(t) \psi(\alpha(t)), \quad t\in
 \sT.
 \end{equation}
Further, a matching pair $(a,b)$ is called Fredholm, if the Toeplitz
operators $T(c)$ and $T(d)$ are Fredholm.

Henceforth the operator $aI$ of multiplication by the function $a\in
L^\infty$ is denoted simply by $a$. Note that if $(a,b)$ is a
matching pair, then the corresponding matrix--function $V(a,b)$
takes the form
  \begin{equation*}
V(a,b)=\left(%
\begin{array}{cc}
 0  & d \\
  -c  & a_{\alpha}^{-1} \\
   \end{array}%
\right).
\end{equation*}
where $(c,d)$ is the subordinated pair for the pair $(a,b)$. In
addition, by \cite{DS:2014} the operator $T(V(a,b))$ can be
represented as the product of three matrix operators, namely,
 \begin{align}
T(V(a,b))& =\left(%
\begin{array}{cc}
 0  & T(d) \\
  -T(c)  & T(a_{\alpha}^{-1}) \\
   \end{array}%
\right)   \nn\\[1ex]
 &= \left(%
\begin{array}{cc}
 -T(d)  & 0 \\
 0  & I \\
   \end{array}%
\right)
\left(%
\begin{array}{cc}
 0  & -I \\
  I  & T(a_{\alpha}^{-1}) \\
   \end{array}%
\right)
\left(%
\begin{array}{cc}
 -T(c)  &  0\\
  0  & I \\
   \end{array}%
\right),  \label{eq16}
\end{align}
where the operator
 \begin{equation*}
   D:=\left(%
\begin{array}{cc}
 0  & -I \\
  I  & T(a_{\alpha}^{-1}) \\
   \end{array}%
\right)
\end{equation*}
in the right--hand side of \eqref{eq16} is invertible and
 \begin{equation*}
   D^{-1} =\left(%
\begin{array}{cc}
 T(a_{\alpha}^{-1})  & I \\
 - I  & 0 \\
   \end{array}%
\right).
\end{equation*}
Further, let us also mention a useful representation for the kernel
of the block Toeplitz operator $T(V(a,b))$ established recently in
the context of classical Toeplitz plus Hankel operators. Thus the
following result holds.

 \begin{prop}[{\cite{DS:2014}}]\label{p3.4}
Let $(a,b)\in L^\infty\times L^\infty$ be a matching pair such that
the operator $T(c)$, $c=ab^{-1}$, is invertible from the right. Then
\begin{equation*}
\ker T(V(a,b))=\Omega(c) \dotplus \widehat{\Omega}(d)
\end{equation*}
 where
 \begin{align*}%\label{eq}
\Omega(c) &:=\left \{   (\vp,0)^T:\vp\in \ker T(c)\right\},\\
\widehat{\Omega}(d) &:=\left \{ (T_r^{-1}(c)
T(a_{\alpha}^{-1})s,s)^T:s\in \ker T(d) \right \},
\end{align*}
and $T_r^{-1}(c)$ is one of the right inverses for the operator
$T(c)$.
 \end{prop}

Thus the inclusion $\vp \in \ker T(c)$ implies that $(\vp,0)^T \in
\ker T(V(a,b))$ and by Lemma \ref{l1} one obtains
  \begin{equation}\label{eq17}
  \begin{aligned}
 \vp-J_{\alpha}QcP\vp \in \ker(T(a)+H_{\alpha}(b)),\\
% \intertext{and}
 \vp+J_{\alpha}QcP\vp\in \ker(T(a)-H_{\alpha}(b)).
\end{aligned}
 \end{equation}
It is even more remarkable that the functions $
\vp-J_{\alpha}QcP\vp$ and $\vp+J_{\alpha}QcP\vp$ belong to the
kernel of the operator $T(c)$ as well.
 \begin{prop}\label{p1}
Assume that $g\in L^\infty$ is a matching function and $f\in \ker
T(g)$. Then $J_{\alpha}QgPf\in \ker T(g)$ and
$(J_{\alpha}QgP)^2f=f$.
  \end{prop}

\textbf{Proof.}
 If $g g_{\alpha}=1$ and $f\in \ker T(g)$, then
  \begin{align*}
  T(g)(J_{\alpha}QgPf)&=PgPJ_{\alpha}QgPf
  =J_{\alpha}Qg_{\alpha}QgPf\\
  &=J_{\alpha}Qg_{\alpha}gPf-J_{\alpha}Qg_{\alpha}PgPf=0,
\end{align*}
i.e. $J_{\alpha}QgPf\in \ker T(g)$. On the other hand, for any $f\in
\ker T(g)$ one has
\begin{align*}%\label{eq}
(J_{\alpha}QgP)^2f & =J_{\alpha}QgPJ_{\alpha}QgP f= Pg_{\alpha}QgPf&\\
&=Pg_{\alpha}gPf- Pg_{\alpha}PgPf =f-Pg_{\alpha}T(g)f=f,
\end{align*}
and we are done. \rbx

Let us introduce the operator
$\mathbf{P}_{\alpha}(g):=J_{\alpha}QgP\left |_{\ker T(g)}\right .$.
Proposition~\ref{p1} implies that $\mathbf{P}_{\alpha}(g):\ker
T(g)\to \ker T(g)$ and $\mathbf{P}_{\alpha}^2(g)=I$. Thus on the
space $\ker T(g)$ the operators
$\mathbf{P}_{\alpha}^{\pm}(g):=(1/2)(I\pm \mathbf{P}_{\alpha}(g))$
are complementary projections, so they generate a decomposition of
$\ker T(g)$. Moreover, using \eqref{eq17} one can formulate the
following result.

 \begin{cor}\label{c1}
Let $(c,d)$ be the subordinated pair for the matching pair $(a,b)\in
L^\infty\times L^\infty$. Then
 \begin{equation*}%\label{eq}
   \ker T(c)=\im \mathbf{P}_{\alpha}^{-}(c)\dotplus\im\mathbf{P}_{\alpha}^{+}(c),\\
   \end{equation*}
   and
  \begin{equation}\label{eq18}
  \begin{aligned}
  & \im \mathbf{P}_{\alpha}^{-}(c)\subset \ker(T(a) +H_{\alpha}(b)) ,\\
   & \im \mathbf{P}_{\alpha}^{+}(c)\subset \ker(T(a)-H_{\alpha}(b)),
\end{aligned}
 \end{equation}
hold.
 \end{cor}
Relations \eqref{eq18} show the influence of the operator $T(c)$ on
the kernels of the operators $T(a)+H_{\alpha}(b)$ and
$T(a)-H_{\alpha}(b)$. Let us now clarify the role of the other
operator--viz. the operator $T(d)$, in the description of the
kernels of the operators $T(a)+ H_{\alpha}(b)$ and $T(a)-
H_{\alpha}(b)$. Assume additionally that the operator $T(c)$ is
invertible from the right. If $s\in \ker T(d)$, then the element
$(T_r^{-1}(c) T(a_{\alpha}^{-1})s,s)^T \in \ker T(V(a,b))$. By Lemma
\ref{l1}, the element
 \begin{equation}\label{eqphi}
  2 \vp_{\alpha}^{\pm}(s):=T_r^{-1}(c)T(a_{\alpha}^{-1}) s \mp J_{\alpha}QcP
   T_r^{-1}(c)T(a_{\alpha}^{-1})s \pm J_{\alpha}Q a_{\alpha}^{-1} s
\end{equation}
belongs to the null space  $\ker(T(a)\pm H_{\alpha}(b))$ of the
corresponding operator $T(a)\pm H_{\alpha}(b)$.

 \begin{lem}\label{l2}
The mapping $s\mapsto \vp_{\alpha}^{\pm}(s)$ is a one-to-one
function from the space $\im \mathbf{P}_{\alpha}^{\pm}(d)$  to the
space $\ker(T(a)\pm H_{\alpha}(b))$.
 \end{lem}
 \textbf{Proof.}
Assuming that $s$ belongs to the kernel of the operator $T(d)$, one
can show that the operator $(1/2)(Pb_{\alpha}P+P
a_{\alpha}J_{\alpha}P)$ sends the element $\vp_{\alpha}^{+}(s)$ into
$\mathbf{P}_{\alpha}^{+}(d)s$ and the operator $(1/2)(Pb_{\alpha}P-P
a_{\alpha}J_{\alpha}P)$ sends the element $\vp_{\alpha}^{-}(s)$ into
$\mathbf{P}_{\alpha}^{-}(d)s$. The proof of these facts is based on
the relations \eqref{eq9} and proceeds similarly to the proof of
\cite[Lemma 3.6]{DS:2014}. Then Lemma \ref{l2} follows.
  \rbx

  \begin{prop}\label{p2}
 Let $(c,d)$ be the subordinated pair for a matching pair $(a,b)\in L^\infty \times L^\infty$.
 If the operator $T(c)$ is right-invertible, then
 \begin{equation}\label{direct}
\begin{aligned}%\label{eq}
   \ker(T(a)+H_{\alpha}(b))& =\vp_{\alpha}^{+}(\im \mathbf{P}_{\alpha}^+(d)) \dotplus\im \mathbf{P}_{\alpha}^-(c), \\
\ker(T(a)-H_{\alpha}(b))& =\vp_{\alpha}^{-}(\im
\mathbf{P}_{\alpha}^-(d)) \dotplus\im \mathbf{P}_{\alpha}^+(c).
\end{aligned}
\end{equation}
Moreover,
 $$
\vp_\alpha^+(\im \mathbf{P}_\alpha^-(d))\subset \im
\mathbf{P}_\alpha^-(c), \quad \vp_\alpha^-(\im
\mathbf{P}_\alpha^+(d))\subset \im \mathbf{P}_\alpha^+(c).
 $$
 \end{prop}

The proof of this proposition is similar to the proof of the
corresponding results of \cite{DS:2014} for classical Toepltz plus
Hankel operators.

 \section{Some classes of generalized Toeplitz plus Hankel operators and
 Coburn--Simonenko theorem.\label{s3}}
The aim of this section is to show how results of the previous
section can be exploited in order to study generalized Toeplitz plus
Hankel operators with generating functions $a$ and $b$ connected in
a special way. Similar classical Toeplitz plus Hankel operators have
been studied in \cite{BE2004, BE2006, Ehr:2004h}. Nevertheless our
approach seems to be more simple and allows us to treat the
operators concerned from a unified point of view. In particular, as
it will be shown below, the operators in question satisfy the
Coburn--Simonenko theorem. Recall that this theorem states that if
$a$ is a non-zero function, then either kernel or cokernel of a
scalar Toeplitz operator $T(g)$, $g\in L^\infty$ is trivial. In
general, generalized Toeplitz plus Hankel operators do not possess
such a property. However, for some classes of operators
$T(a)+H_{\alpha}(b)$ a version of Coburn--Simonenko theorem still
holds. It is worth noting that the approach here does not use
factorization technique.

Let us write the operator $J_{\alpha}: L^p\mapsto L^p$ as follows
 $$
J_{\alpha} \vp(t):= \chi^{-1}(t) \vp (\alpha(t)),
 $$
where
 $$
\chi(t) = \frac{t}{\alpha_-(t)}\,,\quad t\in\sT.
 $$
The function $\chi \in H^\infty$ and has a number of remarkable
properties. Some of them are listed in the lemma below.

 \begin{lem}\label{l3}
 Let the shift $\alpha$ be as above. Then
 \begin{enumerate}
 \item The function $\chi$ is a matching function, i.e.
$\chi \chi_{\alpha}=1$, and $\wind \chi =1$, where $\wind \chi$
denotes the winding number of the function $\chi$.

 \item The function $\chi_\alpha \in \overline{H^\infty}$ and
$\chi_\alpha(\infty)=0$.

\item If $a,b\in L^\infty$ and $n$ is a positive integer, then
   \begin{equation}\label{eqA}
 T(a)+H_{\alpha}(b)=(T(a \chi^{-n})+H_{\alpha}(b \chi^n)) T(\chi^n).
 \end{equation}
 \end{enumerate}
 \end{lem}

\textbf{Proof.} The proof of the first item is a matter of
straightforward computation, whereas the inclusion $\chi_\alpha \in
\overline{H^\infty}$ follows from \eqref{eq7}. The identity
\eqref{eqA} is a consequence of relations \eqref{eq10} and the fact
that $H_\alpha(\chi^n)T(\chi^n)=0$ for any $n\in\sN$.
 \rbx

Now we can establish the main result of this section.

 \begin{thm}\label{t2.1}
Let $a\in GL^\infty$, and let $A$ be any of the operators
$T(a)-H_{\alpha}(a \chi^{-1})$, $T(a)+H_{\alpha}(a \chi )$,
$T(a)+H_{\alpha}(a)$, $T(a)-H_{\alpha}(a)$. Then $\ker A=\{0\}$ or
$\coker A=\{0\}$.
 \end{thm}

 \textbf{Proof.}
Consider first the operator $T(a)+H_{\alpha}(a\chi)$. The duo $(a,a
\chi)$ is a matching pair with the subordinated pair $(c, d)$, where
$c=\chi^{-1}$ and $d=a a_{\alpha}^{-1}\chi$.  Note that $\wind c=-1$
so that the operator $T(c)=T(\chi^{-1})$ is invertible from the
right and $\dim\ker T(\chi^{-1})=1$. Moreover, it is easily seen
that the constant function $\mathbf{e}:=\mathbf{e}(t)=1$, $t\in\sT$,
belongs both to $\ker T(\chi^{-1})$ and $\ker (T(a)-H_{\alpha}(a
\chi))$. Thus $\ker T(\chi^{-1})=\{\lambda \mathbf{e}\}=\im
\mathbf{P}_{\alpha}^-(\chi^{-1})\dotplus \im
\mathbf{P}_{\alpha}^+(\chi^{-1})$, $\lambda\in \sC$. Besides
 \begin{align*}%\label{eq}
&\mathbf{P}_{\alpha}^-(\chi^{-1}) \mathbf{e}  =
(1/2)(\mathbf{e}-J_{\alpha}Q\chi^{-1}P\mathbf{e})\\
&\quad =(1/2)(\mathbf{e}-J_{\alpha}(\chi^{-1}\mathbf{e}))
 =(1/2)(\mathbf{e}-\chi_{\alpha}^{-1}\chi^{-1}\mathbf{e})=0,
\end{align*}
so $\im \mathbf{P}_{\alpha}^-(\chi^{-1})=\{0\}$ and
Proposition~\ref{p2} implies that
\begin{equation}\label{eq3.1.1}
\dim\ker(T(a)+H_{\alpha}(a \chi))=\dim \im \mathbf{P}_{\alpha}^+(d).
\end{equation}
with $d=a a_{\alpha}^{-1}\chi$. Let $\dim \ker T(d)>0$. Then the
Coburn--Simonenko theorem gives  that
\begin{equation*}
 \coker T(\chi^{-1})=\coker T(d)=\{0\}.
 \end{equation*}
Factorization \eqref{eq16} entails that the cokernel of $T(V(a,
a\chi))$ is trivial. By \eqref{eq13} the cokernel of the diagonal
operator $\diag (T(a)+H_{\alpha}(a\chi),T(a)-H_{\alpha}(a \chi))$ is
also trivial, and hence so is the cokernel of $T(a)+H_{\alpha}(a
\chi)$ and that of $T(a)-H_{\alpha}(a \chi)$. On the other hand, if
$\dim\ker T(d)=0$, then \eqref{eq3.1.1} leads to the conclusion that
 \begin{equation*}
\ker (T(a)+H_{\alpha}(a\chi))=\{0\},
 \end{equation*}
so the operator $T(a)+H_{\alpha}(a\chi)$ is subject to the
Coburn--Simonenko theorem.

Consider now the operator $T(a)-H_{\alpha}(\chi^{-1}a)$. The
representation \eqref{eqA} implies that
\begin{equation}\label{eqB}
T(a)-H_{\alpha}(a \chi^{-1})= (T(a \chi^{-1})-H_{\alpha}((a
\chi^{-1})\chi))\cdot T(\chi).
\end{equation}
Setting $b:=a\chi^{-1}$, one rewrites the first operator in the
right-hand side of \eqref{eqB} as
\begin{equation*}
T(a \chi^{-1})-H_{\alpha}((a \chi^{-1})\chi)=T(b)- H_{\alpha}(b
\chi).
\end{equation*}
The operators of the form $T(b)- H_{\alpha}(b \chi)$ have been
already considered, and it was mentioned that the element
$\mathbf{e}$ belongs to the kernel of the operator $T(b)-
H_{\alpha}(b\chi)$. Further, it is not difficult to see that
$\mathbf{e}\notin \im T(\chi)=\im T(\overline{\beta}t-1)$. Indeed,
if $\mathbf{e}\in \im T(\overline{\beta}t-1)$, then there is a
function $h\in H^p$ such that
  \begin{equation}\label{eqC}
1=(\overline{\beta}t-1)h.
 \end{equation}
However, the function $\overline{\beta}t-1$ vanishes at the point
$t=1/\overline{\beta}\in D^+$. The function $h$ admits an analytical
extension into domain $D^+$, whereas relation \eqref{eqC} is still
valid for the domain $D^+$. But this is a contradiction.  Write $d=b
b_\alpha^{-1}\chi$ and assume that $\ker T(d)=\{0\}$. Since
$\mathbf{e}\notin \im T(\chi)$, relation \eqref{eqB} implies that
\begin{equation*}
\ker(T(a)-H_{\alpha}(a \chi^{-1}))=\{0\}.
\end{equation*}
On the other hand, if $\dim\ker T(d)>0$, then
\begin{equation*}
\ker(T(b)-H_{\alpha}((b \chi))=\vp_{\alpha}^-(\im
\mathbf{P}_{\alpha}^-(d)) +\Span\{\mathbf{e}\}
\end{equation*}
and
 $$
 \coker (T(b)-H_{\alpha}(b \chi))=\{0\}.
 $$
Since $\ind T(\overline{\beta}t-1)=-1$, we have that $\im
T(\overline{\beta}t-1) \dotplus \{\lambda \mathbf{e}\}=H^p$,
$\lambda\in \sC$. Let $K$ denote the projection which maps the space
$H^p$ onto $\im T(\overline{\beta}t-1)$ in parallel to $\{\lambda
\mathbf{e}\}$. Each element $s\in \ker (T(b)-H_{\alpha}(b\chi))$ can
be represented as
 $$
s= Ks + (I-K)s,
$$
and $(I-K)s$ is a constant. Thus $(I-K)s \in \ker
(T(b)-H_{\alpha}(b\chi))$ and this implies that $Ks \in \ker
(T(b)-H_{\alpha}(b\chi))$. Using \eqref{eqB}, we get that
 $$
\coker(T(a)-H_{\alpha}(a\chi^{-1}))=0.
$$
The remaining operators $T(a)+H_{\alpha}(a)$ and
$T(a)-H_{\alpha}(a)$ can be considered analogously.
 \rbx

Before concluding this section, let us mention a certain duality for
the operators in Theorem \ref{t2.1}. Along with $\chi(t):=
t/\alpha_-(t)$ consider the function $\Psi(t)= t/\alpha_+(t)$. Then
the following result is true.

 \begin{cor}\label{c3.3}
Let $a\in GL^\infty$, and let $A$ be any of the operators
$T(a)-H_{\alpha}(a_{\alpha} \Psi^{-1})$, $T(a)+H_{\alpha}(a_{\alpha}
\Psi )$, $T(a)+H_{\alpha}(a_{\alpha})$,
$T(a)-H_{\alpha}(a_{\alpha})$. Then $\ker A=0$ or $\coker A=0$.
 \end{cor}

The proof of this corollary directly follows from Theorem \ref{t2.1}
and the representation \eqref{adjoint} for the adjoint of
generalized Toeplitz plus Hankel operator.

\section{A kernel decomposition for Toeplitz operators
 with $\alpha$-matching generating functions.\label{s4}}
In this section we study the kernels of the Toeplitz operators
$T(g)$ in the case where the generating function $g\in L^\infty$ is
an $\alpha$-matching function, i.e. if it satisfies the relation $g
g_{\alpha}=1$. Thus we describe bases in the spaces $\im
\mathbf{P}_\alpha^{\pm}(g)$. These results lay the foundation for
the basis construction of the kernel of generalized Toeplitz plus
Hankel operators with generating matching functions.

Let us start by recalling some properties of Toeplitz operators. It
is well known that Fredholmness of the operator $T(a)$ is closely
connected to Wiener--Hopf factorization of the corresponding
generating function $a$. Let $p>1$, $q>1$ be real numbers such that
$p^{-1}+q^{-1}=1$.
 \begin{defn}
A function $g\in L^\infty$ admits a week Wiener--Hopf factorization
in $H^p$, if it can be represented in the form
 \begin{equation}\label{fac}
    g=g_- t^n g_+,\quad g_-(\infty)=1,
\end{equation}
where $n\in \sZ$, $g_+\in H^q$, $g_+^{-1}\in H^p$, $g_-\in
\overline{H^p}$, $g_-^{-1}\in \overline{H^q}$.
 \end{defn}
The weak Wiener--Hopf factorization of a function $g$ is unique, if
it exists. The functions $g_-$ and $g_+$ are called the
factorization factors, and the number $n$ is the factorization
index. If $g\in L^\infty$ and the operator $T(g)$ is Fredholm, the
function $g$ admits the weak Wiener--Hopf factorization with an
index $n=-\ind T(g)$ \cite{BS:2006, LS1987}. Besides, in this case,
the factorization factors possess an additional property--viz. the
linear operator $g_+^{-1}Pg_-^{-1}I$ defined on $\Span
\{t^k:k\in\sZ_+\}$ can be boundedly extended on the whole space
$H^p$. In the following, such a kind of weak Wiener--Hopf
factorization in $H^p$ is called simply Wiener--Hopf factorization
in $H^p$.

A Wiener--Hopf factorization is called \emph{bounded} if the
factorization factors $g_+,g_+^{-1}$ and $g_-,g_-^{-1}$ belong to
$H^\infty$ and $\overline{H^\infty}$, respectively. Obviously, a
bounded weak Wiener--Hopf factorization is automatically a
Wiener--Hopf factorization which does not depend on $p$. However, it
is worth mentioning that there are continuous non-vanishing on $\sT$
functions which do not admit bounded Wiener--Hopf factorization.
Moreover, let us also recall that if $h_1\in H^\infty$, $h_2\in
\overline{H^\infty}$ and $g\in L^\infty$, then $T(h_2 g
h_1)=T(h_2)T(g)T(h_1)$. The last relation is an immediate
consequence of the already mentioned formula
$T(ab)=T(a)T(b)+H(a)H(\widetilde{b})$, where
$\widetilde{b}:=b(1/t)$.

 \begin{thm}[{\rm see} {\cite[Section 5.5]{BS:2006}}]\label{t1}
If $g\in L^\infty$, then the Toeplitz operator $T(g):H^p\to H^p$,
$1<p<\infty$ is Fredholm and $\ind T(g)=-n$ if and only if the
generating function $g$ admits the Wiener--Hopf factorization
\eqref{fac} in $H^p$.
 \end{thm}
Let us emphasize that Fredholmness of a Toeplitz operator depends on
the space where this operator acts (see \cite{BS:2006, LS1987}) and
in many cases there are efficient formulas to compute the index of
the operator $T(g)$ and, therefore, its defect numbers $\dim\ker
T(g)$ and $\dim\coker T(g)$ due to the Coburn--Simonenko Theorem. We
also recall that one-sided inverses of a Fredholm scalar Toeplitz
operator $T(g)$ can be effectively derived. Thus if the
factorization index $n$ of the function $g$ is non-negative, then
$T(g)$ is left--invertible and the operator $T(t^{-n})T^{-1}(g_0)$,
where $g_0:=a t^{-n}$, is one of the left--inverses for $T(g)$. On
the other hand, if $n\leq0$, then $T(g)$ is right--invertible. For
the sake of convenience, in this paper the notation $T_r^{-1}(g)$
always means the operator $T^{-1}(a_0)T(t^{-n})$, which is one of
right inverses for the operator $T(g)$. Besides, for $n>0$ the
kernel of the operator $T(t^{-n})$ is the linear span of the
monomials $1,t, \cdots, t^{n-1}$, i.e. $ \ker T(t^{-n})=\Span
\{1,t,\cdots, t^{n-1}\}$. Moreover, if $T(g)$ is right--invertible
and $\dim\ker T(g)=\infty$, then $T_r^{-1}(g)$ denotes one of right
inverses of $T(g)$.

Assume now that $g\in L^\infty$ satisfies the matching condition
\eqref{eq15} and the corresponding operator $T(g):H^p\mapsto H^p$ is
Fredholm. We already know that if $\ind T(g)>0$, then
 $$
 \ker T(g)= \im \mathbf{P}_{\alpha}^{+}(g)\dotplus
\im\mathbf{P}_{\alpha}^{-}(g).
$$
Now we want to derive an explicit description of the spaces $\im
\mathbf{P}_{\alpha}^{+}(g)$ and $\im\mathbf{P}_{\alpha}^{-}(g)$, and
the result below is the first step towards this goal.

 \begin{thm}\label{t3}
Assume that $g\in L^\infty$ satisfies the matching condition $g
g_{\alpha}=1$ and the operator $T(g):H^p\mapsto H^p$ is Fredholm.
If\/ $\ind T(g)=n$, $n\in \sZ$, then $g$ can be represented in the
form
\begin{equation}\label{eqPM}
g=\xi  g_{+}\chi^{-n}(g_{+}^{-1})_{\alpha}\,,
\end{equation}
where $g_+$ and $n$ occur in the Wiener--Hopf factorization
 \begin{equation}\label{eqPM1}
g= g_-t^{-n} g_+, \quad g_-(\infty)=1,
 \end{equation}
of the function $g$, whereas $\xi \in \{-1,1\}$ and is defined by
  \begin{equation}\label{eq19}
  \xi= \left ( \frac{\lambda}{\overline{\beta}}\right )^n g_+^{-1}\left ( \frac{1}{\overline{\beta}}\right
).
 \end{equation}
  \end{thm}
\textbf{Proof.} It is clear that $g^{-1}$ possesses a weak
Wiener--Hopf factorization
  $$
g^{-1}=g_-^{-1} t^n g_+^{-1}, \quad g_-^{-1}(\infty)=1
  $$
in $H^q$. But $g^{-1}=g_\alpha$, so that
 \begin{equation}\label{eq20}
 \begin{aligned}%\label{eq}
g_-^{-1} t^n g_+^{-1}&=(g_-)_\alpha \, t^{-n}_\alpha \,
(g_+)_\alpha=(g_-)_\alpha \,(\alpha)^{-n}\, (g_+)_\alpha\\
 &=(g_-)_\alpha\,
\alpha_-^{-n}\,t^n \,\alpha_+^{-n}\, (g_+)_\alpha = \eta_- t^n
\eta_+,
\end{aligned}
 \end{equation}
where
  $$
\eta_+=(g_-)_\alpha \, \alpha_+^{-n}, \quad
\eta_-=(g_+)_\alpha\,\alpha_-^{-n}.
  $$
It follows from \eqref{eq5}, \eqref{eq7} and from the properties of
the factorization factors $g_-$ and $g_+$ that $\eta_+\in H^p$,
$\eta_+^{-1}\in H^q$, $\eta_-\in \overline{H^q}$ and $\eta_-^{-1}\in
\overline{H^p}$. Therefore, comparing representations \eqref{eqPM1}
and \eqref{eq20}, one obtains
 $$
 g_-^{-1} \,\eta_-^{-1}=  \eta_+ g_+=\xi
 $$
where $\xi$ is a complex number. Indeed, the product in the
left-hand side belongs to $\overline{H^1}$ whereas the right-hand
side is in $H^1$, and it is well-known that $\overline{H^1}\cap
H^1=\sC$. Thus
 \begin{equation*}
g_-^{-1}=\xi \, \eta_-, \quad g_+=\xi \eta_+^{-1},
\end{equation*}
so that
 $$
g=\xi^{-1}\eta_-^{-1} t^{-n} g_+=\xi^{-1}\,(g_+)_\alpha^{-1}\,
\alpha_+^{-n}\, t^{-n}\, g_+ =\xi^{-1}(g_+)_\alpha^{-1}\,
\chi^{-n}\, g_+.
 $$
Using this identity one can represent $g_\alpha$ in the form
$g_\alpha=\xi^{-1}\,g_+^{-1} \chi^n\,(g_+)_\alpha$, which leads to
the equation $1=g g_\alpha=(\xi^{-1})^2$, i.e $\xi=-1$ or $\xi=1$.
On the other hand, from $g_-=\xi^{-1}\eta_-^{-1}$ and
$g_-(\infty)=1$ we get
 $$
1=\frac{1}{\xi} \lim_{t\to\infty}
(g_+)_\alpha^{-1}(t)(\alpha_-(t))^n=\frac{1}{\xi} \,  g_+^{-1}\left
(\frac{1}{\overline{\beta}}\right )  \left
(\frac{\lambda}{\overline{\beta}} \right )^n,
 $$
and the relation \eqref{eq19} follows.
 \rbx

 \begin{rem}
Factorization \eqref{eqPM} has been mentioned in \cite{KLR:2009}
without proof but the factor $\xi$ is missing there.
 \end{rem}
Now we can introduce the following definition.
 \begin{defn}
The number $\xi$ defined by the relation \eqref{eq19} is called the
$\alpha$-factorization signature, or simply, $\alpha$-signature of
$g$ and is denoted by $\boldsymbol\sigma_{\!\alpha}(g)$.
 \end{defn}
The $\alpha$-signature plays an important role in the construction
of bases of the kernels of generalized Toeplitz plus Hankel
operator. In order to determine it, one has to evaluate the
corresponding factorization factor at a point of the complex plane
$\sC$. In general situation, this is not an easy task at all.
Nevertheless, for some classes of generating functions $g$, this
specific characteristic can be easily found and such a possibility
is discussed later on.

  \begin{thm}\label{t2}
Let $g\in L^\infty$ be an $\alpha$-matching function such that the
operator $T(g):H^p\to H^p$ is Fredholm and $n:=\ind T(g)>0$. If
$g=g_-t^{-n} g_+$, $g_-(\infty)=1$ is the corresponding Wiener--Hopf
factorization of $g$ in $H^p$, then  the following systems of
functions $\cB_{\alpha}^{\pm}(g)$ form bases in the spaces $\im
\mathbf{P}_{\alpha}^{\pm}(g)$:
  \begin{enumerate}
    \item If $n=2m$, $m\in \sN$, then
    \begin{equation*}
\cB_{\alpha}^{\pm}(g):= \{g_+^{-1}(\chi^{m-k-1}\pm
\boldsymbol\sigma_{\!\alpha} (g)\chi^{m+k}):k=0,1,\cdots, m-1
 \},
    \end{equation*}
     and
     \begin{equation*}
      \dim\im \mathbf{P}_{\alpha}^{\pm}(g)=m.
     \end{equation*}
    \item If $n=2m+1$, $m\in \sZ_+$, then
    \begin{equation*}
\cB_{\alpha}^{\pm}(g):= \{g_+^{-1}(\chi^{m+k}\pm
\boldsymbol\sigma_{\!\alpha} (g)\chi^{m-k}):k=0,1,\cdots, m
 \}\setminus \{0\},
    \end{equation*}
    \begin{equation*}
     \dim\im \mathbf{P}_{\alpha}^{\pm}(g)=m +\frac{1\pm\boldsymbol\sigma_{\!\alpha}(g)}{2}
     \,.
    \end{equation*}
   \end{enumerate}
\end{thm}
  \begin{rem}
If $n=2m+1$, then the zero element belongs to one of the sets
$\cB_{\alpha}^{+}(g)$ or $\cB_{\alpha}^{-}(g)$. Namely, for $k=0$
one of the terms $\chi^m(1\pm\boldsymbol\sigma_{\!\alpha}(g))$ is
equal to zero.
  \end{rem}
\textbf{Proof.} [Proof of Theorem \ref{t2}] Observe that the
restriction of the operators $Pg_-I$
 and $Pg_-^{-1}I$ on $\ker T(t^{-n})=\Span \{\chi^0, \chi,\cdots,
 \chi^{n-1} \}$ map $\ker T(t^{-n})$ into $\ker T(t^{-n})$, and on the
 space $\ker T(t^{-n})$ the above operators are inverses to each
 other.

Thus the elements $s_j=Pg_-\chi^j$, $j=0,1,\cdots, n-1$ belong to
$\ker T(t^{-n})$ and
\begin{equation*}
T^{-1}(g_0)s_j=g_+^{-1} Pg_-^{-1}s_j=g_+^{-1}\chi^j,
\end{equation*}
where $g_0:=gt^n$. Here we used the relation
$T^{-1}(g_0)=g_+^{-1}Pg_-^{-1}$ and the fact that $T^{-1}(g_0)
s_j\in\ker T(g)$.  Moreover, the set
$\{T^{-1}(g_0)s_j:j=0,\cdots,n-1 \}$ constitutes a basis in $\ker
T(g)$. Now one can use the factorization \eqref{fac} and write
 \begin{align*}%\label{eq}
&\mathbf{P}_{\alpha}(g) (T^{-1}(g_0) s_j)
=J_{\alpha}QgPT^{-1}(g_0)s_j\\
& \phantom{abc}=
J_{\alpha}Q\boldsymbol\sigma_{\!\alpha}(g)(g_+^{-1})_\alpha
\chi^{-n}g_+P g_+^{-1}\chi^j=\boldsymbol\sigma_{\!\alpha}(g)
g_+^{-1} \chi^{n-j-1},
\end{align*}
which leads to the representation
\begin{equation*}
\mathbf{P}_{\alpha}^{\pm}(g)( T^{-1}(g_0) s_j)=
\frac{1}{2}g_+^{-1}(\chi^j\pm
\boldsymbol\sigma_{\!\alpha}(g)\chi^{n-j-1}), \quad
j=0,1,,\cdots,n-1.
\end{equation*}
Assume now that $n=2m$, $m\in\sN$. If $j\in\{0,1,\cdots,m-1\}$, it
can be rewritten as $j=m-k-1$ with some $k\in\{0,1,\cdots,m-1\}$ and
vice versa. Hence
\begin{equation*}
\chi^j\pm \boldsymbol\sigma_{\!\alpha}(g)\chi^{n-j-1} =
\chi^{m-k-1}\pm \boldsymbol\sigma_{\!\alpha}(g)\chi^{m+k}, \quad
j=0,1,\cdots,m-1.
 \end{equation*}
On the other hand, if $j\geq m$, then $j=m+k$ for a
$k\in\{0,1,\cdots, m-1\}$, and
 $$
\chi^j\pm \boldsymbol\sigma_{\!\alpha}(g)\chi^{n-j-1} =
\chi^{m+k}\pm \boldsymbol\sigma_{\!\alpha}(g)\chi^{m-k-1}=
\pm\boldsymbol\sigma_{\!\alpha}(g) (\chi^{m-k-1}\pm
\boldsymbol\sigma_{\!\alpha}(g)\chi^{m+k}).
 $$
Thus one obtains, maybe up to the factor $-1$,  the same system of
function that is derived for $j\in\{0,1,\cdots, m-1\}$. So we
conclude that if $n=2m$, then
 \begin{equation*}
\dim \im \mathbf{P}_{\alpha}^{\pm}(g)=m,
 \end{equation*}
and assertion (i) is shown.

The proof of assertion (ii) is similar to that of assertion (i).
 \rbx

As was already mentioned, the evaluation of $\alpha$-signature is a
difficult problem. However,  there are functions $g\in L^\infty$,
the $\alpha$-signature of which can be easily determined. Consider,
for example, a function $g\in L^\infty$ continuous at one of the
fixed points $t_{\pm} =(1\pm \lambda)/\overline{\beta}\in\sT$ of the
mapping $\alpha$. At any fixed point, such a function $g$ can take
only one of the two values, namely, $+1$ or  $-1$, which leads to
the following results.

  \begin{prop}\label{p3}
Let $g\in L^\infty$ be an $\alpha$-matching function such that
\begin{enumerate}
    \item The operator $T(g):H^p\to H^p$ is Fredholm with index $n$.
    \item The function $g$ is continuous at the point $t_+$ or $t_-$.
\end{enumerate}
Then
 $$\boldsymbol\sigma_{\!\alpha}(g)=g(t_+) \, \text{ or } \,
\boldsymbol\sigma_{\!\alpha}(g)=g(t_-)(-1)^n.
 $$
  \end{prop}

\textbf{Proof.}
 Assume for definiteness that the function $g$ is continuous at the
 point $t_+$. Condition (i) ensures that $g$ admits a
Wiener--Hopf factorization in $H^p$,
\begin{equation*}
g=g_- t^{-n} g_+, \quad g(\infty)=1,
\end{equation*}
and $g(t_+) \in \{ -1,1\}$. Now we approximate the function $g$ as
follows. Chose an $\ve>0$ and an open arc $\sT_{\ve}\subset \sT$
such that
 \begin{enumerate}
    \item The point $t_+$ belongs to the arc $\sT_{\ve}$.
    \item $\alpha (\sT_{\ve})=\sT_{\ve}$ that is that $\alpha$ is a
    homomorphism of $\sT_{\ve}$.
    \item  If $g_{\ve}$ denotes the function
     $$
g_{\ve}(t) := \left \{
 \begin{array}{ll}
g(t_+) & \, \text{ if }\, t\in \sT_{\ve}\\[1ex]
g(t) & \, \text{ otherwise }
 \end{array},
 \right.
 $$
 then
  $$
||g-g_{\ve}||<\ve.
  $$
 \end{enumerate}
It is clear that such an arc $\sT_{\ve}$ exists and that
$g_{\ve}(g_{\ve})_{\alpha}=1$. If $\ve$ is small enough, then the
operator $T(g_{\ve})$ is also Fredholm on $H^p$ with the same index
$n$. In this case, the function $g_{\ve}$ admits a Wiener--Hopf
factorization in $H^p$,
 $$
g_{\ve}=g_{\ve,-}t^{-n} g_{\ve,+}, \quad g_{\ve}(\infty)=1.
 $$

Since $g_{\ve}$ is H\"older continuous in a neighbourhood of the
point $t_+$, Corollary 5.15 of \cite{LS1987} indicates that the
functions $g_{\ve,-}$ and $g_{\ve,+}$ are also H\"older continuous
in a neighbourhood of $t_+$.  By Theorem \ref{t3}, the function
$g_{\ve}$ can be written in the form
\begin{equation}\label{eq5.1}
g_{\ve}=\boldsymbol\sigma_{\!\alpha}(g_{\ve}) g_{\ve,+}
(g_{\ve,+}^{-1})_{\alpha}\, \chi^{-n} .
\end{equation}
But $g_{\ve}(t_+)=\boldsymbol\sigma_{\!\alpha}(g_{\ve})
g_{\ve,+}(t_+) (g_{\ve,+}^{-1})_{\alpha}(t_+)\, \chi^{-n}(t_+) =
\boldsymbol\sigma_{\!\alpha}(g_{\ve})$ because $\alpha(t_+)=t_+$ and
$\chi(t_+)=1$. Moreover, we have
\begin{align*}%\label{eq}
\boldsymbol\sigma_{\!\alpha}(g)  & = \left (
\frac{\lambda}{\overline{\beta}} \right )^n
 g_+^{-1}\left (\frac{1}{\overline{\beta}}\right ) , \\
\boldsymbol\sigma_{\!\alpha}(g_{\ve}) & = \left (
\frac{\lambda}{\overline{\beta}} \right )^n
 g_{\ve,+}^{-1}\left (\frac{1}{\overline{\beta}}\right ).
\end{align*}
Our aim now is to show that if $\ve$ is small enough, then
$\boldsymbol\sigma_{\!\alpha}(g_{\ve})=\boldsymbol\sigma_{\!\alpha}(g)$.
Consider the functions
 $$
gt^n=g_- g_+, \quad g_{\ve}t^n= g_{\ve,-} g_{\ve,+},
 $$
and note that the Toeplitz operators $T(gt^n)$ and $T(g_{\ve}t^n)$
are invertible on $H^p$. Let us also assume that $\ve$ is chosen so
small that
  \begin{equation}\label{EqNorm}
  ||T^{-1}(gt^n)-T^{-1}(g_{\ve}t^n)||< \left |
  \frac{\lambda}{\overline{\beta}}
  \right |^{-n} s_p^{-1},
 \end{equation}
where $s_p$  denotes the norm of the linear bounded functional
$h\mapsto h(1/\overline{\beta})$, $h\in H^p$. Further, consider two
uniquely solvable equations
 $$
T(gt^n)h=1,\quad T(g_{\ve}t^n)k=1.
 $$
Using the above mentioned Wiener--Hopf factorization one obtains
\begin{align*}%\label{eq}
    h&=T^{-1}(gt^n)1= g_{+}^{-1} P g_{-}^{-1}1=g_{+}^{-1},\\
k&=T^{-1}(g_{\ve}t^n)1= g_{\ve,+}^{-1}P
g_{\ve,-}^{-1}1=g_{\ve,+}^{-1},
\end{align*}
and
 $$
||g_{+}^{-1}-g_{\ve,+}^{-1}||\leq
||T^{-1}(gt^n)-T^{-1}(g_{\ve}t^n)||.
 $$
It follows that
\begin{align*}
    |\boldsymbol\sigma_{\!\alpha}(g)
    -\boldsymbol\sigma_{\!\alpha}(g_{\ve})|& = \left |
\left ( \frac{\lambda}{\overline{\beta}}\right )^n  \left ( g_+^{-1}
\left ( \frac{1}{\overline{\beta}} \right )-  g_{\ve,+}^{-1} \left (
\frac{1}{\overline{\beta}}\right )\right ) \right | \\
 &\leq\left | \frac{\lambda}{\overline{\beta}}\right |^n s_p \,\,
||T^{-1}(gt^n)-T^{-1}(g_{\ve}t^n)||,
\end{align*}
and relation \eqref{EqNorm} leads to the inequality
 $$
|\boldsymbol\sigma_{\!\alpha}(g)
    -\boldsymbol\sigma_{\!\alpha}(g_{\ve})| <1.
 $$
Therefore,
 $$
\boldsymbol\sigma_{\!\alpha}(g)
    =\boldsymbol\sigma_{\!\alpha}(g_{\ve}),
 $$
and we have
 $$
\boldsymbol\sigma_{\!\alpha}(g)
    =\boldsymbol\sigma_{\!\alpha}(g_{\ve})= g_{\ve}(t_+)=g(t_+).
 $$
If $g$ is continuous at the point $t_-$, the proof is analogous but
one has take into account that $\chi(t_-)=-1$.
 \rbx

\begin{rem}\label{r4}
In the proof we used the fact that for a fixed $\xi\in D^+$ the
mapping $h\mapsto h(\xi)$ is a bounded linear functional on $H^p$.
The boundedness of this functional can be easily seen. Indeed, by
Caushy's integral formula for any polynomial $P_n$, one has
 $$
P_n(\xi) =\frac{1}{2\pi}\int_0^{2\pi}
\frac{P_n(e^{i\theta})e^{i\theta}}{e^{i\theta}-\xi}\, d\theta.
 $$
Using H\"older inequality, one obtains that on the set of all
polynomials the linear functional $h\mapsto h(\xi)$ is bounded in
the $H^p$-norm. Since this set is dense in $H^p$, our claim follows.
 \end{rem}
Note that the two cases above exhaust all the situations possible.

\section{Bases of the kernels and cokernels of generalized Toeplitz plus
Hankel operators.\label{s5}}
 In this section the structure of the
kernel and cokernel of generalized Toeplitz plus Hankel operator
$T(a)+H_{\alpha}(b)$ is described. The operators in question are
studied under the condition that their generating functions $a, b\in
L^\infty$ constitute a Fredholm matching pair $(a,b)$. Recall that
if a matching pair $(a,b)$ is Fredholm, then it follows from
\eqref{eq13} and \eqref{eq16} that $T(a)+H_{\alpha}(b)$ and
$T(a)-H_{\alpha}(b)$ are Fredholm operators. Set $\kappa_1:=\ind
T(c)$, $\kappa_2:=\ind T(d)$ and let $\sZ_-$ refer to the set of all
negative integers.
 \begin{thm}\label{t5.1}
Assume that  $(a,b)\in L^\infty\times L^\infty$ is a Fredholm
matching pair. Then
 \begin{enumerate}
    \item If $(\kappa_1,\kappa_2)\in \sZ_+\times \sN$, then the
    operators $T(a)+H_{\alpha}(b)$ and $T(a)-H_{\alpha}(b)$ are invertible from the
    right and
    \begin{align*}%\label{eq}
\ker (T(a)+H_{\alpha}(b))&=\im\mathbf{P}_{\alpha}^{-}(c)\dotplus \vp_{\alpha}^+(\im \mathbf{P}_{\alpha}^{+}(d)), \\
\ker (T(a)-H_{\alpha}(b))&=\im\mathbf{P}_{\alpha}^{+}(c)\dotplus
\vp_{\alpha}^-(\im \mathbf{P}_{\alpha}^{-}(d)),
\end{align*}
where the spaces $\im\mathbf{P}_c^{\pm}$ and $\im\mathbf{P}_d^{\pm}$
are described in Theorem \ref{t2}, and the mappings
$\vp_{\alpha}^{\pm}$ are defined by \eqref{eqphi}.

     \item  If $(\kappa_1,\kappa_2)\in  \sZ_-\times (\sZ\setminus \sN)$, then the
    operators $T(a)+H_{\alpha}(b)$ and $T(a)-H_{\alpha}(b)$ are invertible from the
    left and
    \begin{align*}%\label{eq}
\coker (T(a)+H_{\alpha}(b))&=\im
\mathbf{P}_{\alpha}^{-}(\overline{d})\dotplus \vp_{\alpha}^+(\im\mathbf{P}_{\alpha}^{+}(\overline{c})), \\
\coker (T(a)-H_{\alpha}(b))&=\im
\mathbf{P}_{\alpha}^{+}(\overline{d})\dotplus
\vp_{\alpha}^-(\im\mathbf{P}_{\alpha}^{-}(\overline{c})),
\end{align*}
and $\im \mathbf{P}_{\overline{d}}^{\pm}=\{0\}$ for $\kappa_2=0$.

   \item If $(\kappa_1,\kappa_2)\in \sZ_+\times (\sZ\setminus \sN)$, then
    \begin{align*}%\label{eq}
\ker (T(a)+H_{\alpha}(b))=\im\mathbf{P}_{\alpha}^{-}(c),  & \quad  \coker (T(a)+H_{\alpha}(b)=\im \mathbf{P}_{\alpha}^{+}(\overline{d}), \\
\ker (T(a)-H_{\alpha}(b))=\im\mathbf{P}_{\alpha}^{+}(c)  & \quad
\coker (T(a)-H_{\alpha}(b)=\im
\mathbf{P}_{\alpha}^{-}(\overline{d}).
\end{align*}

 \end{enumerate}
 \end{thm}

\textbf{Proof.} Note that all results concerning the kernels of the
operators under consideration follow from Theorem \ref{t2} and
representations \eqref{direct}. In order to describe the cokernels
of the corresponding operators, let us recall that $\coker (T(a)\pm
H_{\alpha}(b)):=\ker (T(a)\pm H_{\alpha}(b))^*$. Moreover, $(T(a)\pm
H_{\alpha}(b))^*=T(\overline{a})\pm
H_{\alpha}(\overline{b}_{\alpha})$ and the duo $(\overline{a},
\overline{b}_{\alpha})$ is a matching pair with the subordinated
pair $(\overline{d}, \overline{c})$. Further, if $g\in L^\infty$ and
$\Ind c=r$, then $\Ind \overline{c}=-r$, and the description of the
cokernel immediately follows from the previous results for the
kernels of Toeplitz plus Hankel operators.
 \rbx

It remains to consider the case  $(\kappa_1,\kappa_2)\in \sZ_-\times
\sN$. This situation is more involved and factorization \eqref{eq16}
already indicates that for $\kappa_2>0$, the dimension of the kernel
$\diag(T(a)+H_{\alpha}(b),T(a)-H_{\alpha}(b))$ may be smaller than
$\kappa_2$. To treat this case, consider a number $n\in\sN$ such
that
\begin{equation*}
0\leq 2n+\kappa_1\leq 1.
\end{equation*}
Such an $n$ is uniquely defined and
\begin{equation*}
2n+\kappa_1 =\left\{%
\begin{array}{ll}
   0, & \hbox{if\;} \kappa_1 \; \hbox{is even,}  \\
    1, &\hbox{if\;} \kappa_1 \; \hbox{is odd.} \\
\end{array}%
\right.
\end{equation*}
Now one can use the relation \eqref{eqA} and represent the operators
$T(a)\pm H_{\alpha}(b)$ in the form
 \begin{equation}\label{eq6.1}
T(a)\pm H_{\alpha}(b)= ( T(a\chi^{-n})\pm
H_{\alpha}(b\chi^n))T(\chi^n).
\end{equation}
But $(a\chi^{-n}, b\chi^{n})$ is a matching pair with the
subordinated pair $(c\chi^{-2n}, d)$. Hence, the operators
$T(a\chi^{-n})\pm H_{\alpha}(b\chi^n)$ are subject to assertion (i)
of Theorem \ref{t5.1}. Thus they are right-invertible, and if
$\kappa_1$ is even, then
 \begin{equation}\label{eq6.2}
\begin{aligned}
    \ker (T(a\chi^{-n})+ H_{\alpha}(b\chi^n) )=\vp_{\alpha}^+(\im \mathbf{P}_{\alpha}^{+}(d)),\\
 \ker (T(a\chi^{-n})- H_{\alpha}(b\chi^n) )=\vp_{\alpha}^-(\im \mathbf{P}_{\alpha}^{-}(d)),
\end{aligned}
\end{equation}
and if $\kappa_1$ is odd, then
\begin{equation}\label{eq6.3}
\begin{aligned}
    \ker (T(a\chi^{-n})+ H_{\alpha}(b\chi^n) )= \frac{1-\boldsymbol\sigma_{\!\alpha}(c)}{2}c_+^{-1}\sC
\dotplus  \vp_{\alpha}^+(\im \mathbf{P}_{\alpha}^{+}(d)),\\
 \ker (T(a\chi^{-n})- H_{\alpha}(b\chi^n) )=\frac{1+\boldsymbol\sigma_{\!\alpha}(c)}{2}c_+^{-1}\sC
 \dotplus \vp_{\alpha}^-(\im \mathbf{P}_{\alpha}^{-}(d)),
\end{aligned}
\end{equation}
where the mappings $\vp_{\alpha}^{\pm}$ depend on the functions
$a\chi^{-n}$ and $b\chi^n$.

 \begin{thm}\label{t4}
Let $(\kappa_1,\kappa_2)\in \sZ_-\times \sN$. Then
\begin{enumerate}
    \item If $\kappa_1$ is odd, then
     \begin{align*}%\label{eq}
 &\ker(T(a)\pm  H_{\alpha}(b)) = \\
  & \quad T(\chi^{-n}) \left (\left\{
 \frac{1\mp\boldsymbol\sigma_{\!\alpha}(c)}{2}c_+^{-1}\sC
 \dotplus \vp_{\alpha}^{\pm}(\im \mathbf{P}_{\alpha}^{\pm}(d))\right\} \cap \im T(\chi^n) \right );%\\
% &= \left\{u\in  \left
%\{\frac{1\pm\boldsymbol\sigma_{\!\alpha}(c)}{2}c_+^{-1}\sC
 %\dotplus \vp_{\alpha}^{\pm}(\im \mathbf{P}_{\alpha}^{\pm}(d))\right \} :
% \widehat{u}_0=\cdots =\widehat{u}_{n-1}=0\right \},
\end{align*}
%where $\widehat{u}_k, k=0,1,\cdots, n-1$ are the corresponding
%Fourier coefficients of the function $u$,
    \item If $\kappa_1$ is even, then
     \begin{align*}%\label{eq}
&\ker(T(a)\pm  H_{\alpha}(b)) = T(\chi^{-n}) \left (
\left\{\vp_{\alpha}^{\pm}(\im \mathbf{P}_{\alpha}^{\pm}(d))
\right\} \cap \im T(\chi^n) \right ),%\\
% &= \left\{ u\in  \vp_{\alpha}^{\pm}(\im \mathbf{P}_{\alpha}^{\pm}(d)) :
 %\widehat{u}_0=\cdots =\widehat{u}_{n-1}=0\right \}.
 \end{align*}
\end{enumerate}
and the mappings $\vp_{\alpha}^{\pm}$ depend on the functions
$a\chi^{-n}$ and $b\chi^n$.
 \end{thm}

\begin{rem}
 Using the fact that the system $\{t^{k-1} \alpha_+^{-n}\}_{k=1}^n$
forms a basis of $\ker T^*(\chi^n)=\ker T(\overline{\chi}^n) =\ker
T(t^{-n}\alpha_+^n)$, it is easily seen that $h\in H^p$ belongs to
$\im T(\chi^n)$ if and only if
 $$
\widehat{(h \alpha_-^n)}_i=0, \quad i=0,1, \cdots, n-1,
 $$
where $\widehat{(h \alpha_-^n)}_i$ denotes the $i$-th Fourier
coefficient of the function $h \alpha_-^n$.
 \end{rem}

\textbf{Proof.} It follows immediately from representations
\eqref{eq6.1}--\eqref{eq6.3}. \rbx

Theorem \ref{t4} can also be used to describe the cokernels of the
operators $T(a)\pm H_{\alpha}(b)$ in the situation where
$(\kappa_1,\kappa_2)\in \sZ_-\times \sN$. Indeed, recalling that
$(T(a)\pm H_{\alpha}(b))^*=T(\overline{a})\pm
H_{\alpha}(\overline{b}_{\alpha})$, and $(\overline{d},
\overline{c})$ is the subordinated pair for
$(\overline{a},\overline{b}_{\alpha})$, one can note that the
operators $T(\overline{d})$ and $T(\overline{c})$ are also Fredholm
and
\begin{equation*}
\ind T(\overline{d})=-\kappa_2, \quad \ind
T(\overline{c})=-\kappa_1,
\end{equation*}
so $(-\kappa_2,\kappa_1)\in \sZ_-\times \sN$. Therefore, Theorem
\ref{t4} applies and we can formulate the following result.

 \begin{thm}\label{t5}
Let $(\kappa_1,\kappa_2)\in \sZ_-\times \sN$, and let $m\in\sN$
satisfy the requirement
\begin{equation*}
1\geq 2m-\kappa_2\geq0.
\end{equation*}
Then
\begin{enumerate}
    \item If $\kappa_2$ is odd, then
     \begin{align*}%\label{eq}
 &\coker(T(a)\pm H_{\alpha}(b)) = \\
& \qquad T(\chi^{-m}) \left (\left\{
 \frac{1\mp\boldsymbol\sigma_{\!\alpha}(\overline{d})}{2}\overline{d_-^{-1}}\sC
 \dotplus \vp_{\alpha}^{\pm}(\im \mathbf{P}_{\alpha}^{\pm}(\overline{c}))\right\} \cap \im T(\chi^m) \right ).
\end{align*}

    \item If $\kappa_2$ is even, then
     \begin{align*}%\label{eq}
 \ker(T(a)\pm H_{\alpha}(b)) &= T(\chi^{-m}) \left (\left\{\vp_{\alpha}^{\pm}
  (\im \mathbf{P}_{\alpha}^{\pm}(\overline{c}))\right\} \cap \im
 T(\chi^m) \right ),
 \end{align*}
and the mappings $\vp_{\alpha}^{\pm}$ depend on
$\overline{a}\chi^{-m}$ and $\overline{b}_{\alpha}\chi^m$.
\end{enumerate}
\end{thm}
Thus Theorems \ref{t5.1}-\ref{t5} offer an explicit description of
the kernels and cokernels of the operators under consideration. On
the other hand, the above approach can be also used to find
generalized Toeplitz plus Hankel operators which are subject to
Coburn--Simonenko theorem. Recall that some operators possessing
this property have been already studied in Section~\ref{s3}.

  \begin{prop}\label{c5}
Let $(a,b)\in L^\infty \times L^\infty$ be a matching pair with the
subordinated pair $(c,d)$, and let $T(c)$ be a Fredholm operator.
Then:
\begin{enumerate}
    \item If\/ $\ind T(c)=1$ and $\boldsymbol\sigma_{\!\alpha}(c)=1$, then $\ker(T(a)+H_{\alpha}(b))=\{0\}$
     or $\coker(T(a)+H_{\alpha}(b))=\{0\}$.
    \item If\/ $\ind T(c)=-1$ and $\boldsymbol\sigma_{\!\alpha}(c)=1$, then $\ker(T(a)-H_{\alpha}(b))=\{0\}$ or
    $\coker(T(a)-H_{\alpha}(b))=\{0\}$.
    \item If\/ $\ind T(c)=0$, then $\ker(T(a)\pm H_{\alpha}(b))=\{0\}$ or $\coker(T(a)\pm
    H_{\alpha}(b))=\{0\}$.
\end{enumerate}
 \end{prop}
\textbf{Proof.} The proof of the last theorem is similar to the
proof of Theorem~\ref{t2.1}, but in the proofs of assertions (i) and
(ii) one, respectively, has to use the fact that $c_+^{-1}\in \ker
(T(a)+H_{\alpha}(b))$ and $c_+^{-1}\in \ker (T(a
\chi^{-1})-H_{\alpha}(b\chi))$. These inclusions can be verified by
straightforward computations. Thus consider, for example, the
expression $(T(a \chi^{-1})-H_{\alpha}(b\chi))c_+^{-1}$. Taking into
account factorization \eqref{eqPM} and using the relations
$c_-=(c_+^{-1})_{\alpha}$ and  $a=bc$, one obtains
  \begin{align*}
(T(a \chi^{-1})-H_{\alpha}(b\chi))c_+^{-1}& = Pbc \chi^{-1}c_+^{-1}-Pb\chi QJ_{\alpha}c_+^{-1}\\
 &=Pbc_+\chi c_- \chi^{-1}c_+^{-1}-Pb\chi(c_+^{-1})_{\alpha}
\chi^{-1}=0,
  \end{align*}
and the inclusion $c_+^{-1}\in \ker (T(a
\chi^{-1})-H_{\alpha}(b\chi)$ is proved.
 \rbx

 \begin{rem}
Let us emphasize that Proposition \ref{c5} is valid without any
assumption about Fredlomness, semi-Fredholness, or even normal
solvability of the operators $T(a)\pm H_{\alpha}(b)$. However, if
one of these operators is Fredholm and its index is known, then
Proposition \ref{c5} allows one to compute the kernel and cokernel
dimension of the operator under consideration. For the case of
piecewise continuous generating function see also Section \ref{s6}.
 \end{rem}

 \begin{cor}\label{c6}
Let $b\in L^\infty$ be a matching function. If\/ $T(b_{\alpha})$ is
a Fredholm operator, then:
\begin{enumerate}
    \item If\/ $\ind T(b_{\alpha})=1$,
    and $\boldsymbol\sigma_{\!\alpha}(b_{\alpha})=1$, then $\ker(I+H_{\alpha}(b))$  or
    $\coker(I+H_{\alpha}(b))$ is trivial.
    \item If\/ $\ind T(b_{\alpha})=-1$,
    and $\boldsymbol\sigma_{\!\alpha}(b_{\alpha})=1$, then $\ker(I-H_{\alpha}(b))$ or $\coker(I-H_{\alpha}(b))$ is trivial.
    \item If\/ $\ind T(b_{\alpha})=0$,
    then $\ker(I\pm H_{\alpha}(b))$ or $\coker(I\pm H_{\alpha}(b))$ is trivial.
\end{enumerate}
 \end{cor}
These results are direct consequences of  Proposition \ref{c5},
since if $b$ is a matching function, then $(1,b)$ is a matching pair
with the subordinate pair $(b_{\alpha},b)$.

In conclusion of this section, we would like to mention that in
certain cases the condition of Fredholmness of the operator $T(d)$
can be dropped. However, we are not going to pursue this matter
here.

% \end{document}

\section{Fredholmness of generalized Toeplitz plus Hankel operators with piecewise
continuous generating functions\label{s6}}

Assume $\alpha$ is the Carleman shift \eqref{eq4} changing the
orientation, i.e. $|\beta|>1$. As was already mentioned, the mapping
$\alpha:\sT\to\sT$ has two fixed points, namely
 $$
t_{\pm}=\frac{1\pm i \sqrt{|\beta|^2-1}}{\overline{\beta}}.
$$
By $\sT^+_\alpha$ and $\sT^-_\alpha$ we denote the closed arcs of
$\sT$ which, respectively, join $t_+$ with $t_-$ and $t_-$ with
$t_+$ and inherit the orientation of $\sT$. Further, let us
introduce the functions
 \begin{align*}
 \nu_p(y)&:=\frac{1}{2}\left ( 1+\coth\left (\pi
\left(y+\frac{i}{p}\right ) \right )\right
    ), \quad %\\[1ex]
 h_p(y):= \sinh^{-1}\left (\pi \left(y+\frac{i}{p}\right )
\right ),
\end{align*}
where $y\in \overline{\sR}$, and $\overline{\sR}$ refers to the
two-point compactification of $\sR$. Set $\tp:=\sT_\alpha^+\setminus
\{ t_+,t_-\}$.

 \begin{thm}\label{tTH}
 If $a,b\in PC$, then
 the operator $T(a)+H_{\alpha}(b)$ is Fredholm if and only if the matrix
  \begin{align*}
    &\smb(T(a)+H_{\alpha}(b))(t,y):= \nn\\[1ex]
  &  \left(\!\!\!%
\begin{array}{c@{\hspace{-1mm}}c}
a(t+0)\nu_p(y)+a(t-0)(1-\nu_p(y))   & \D \frac{b(t+0)-b(t-0)}{2i} \,h_p(y)\\
\D \frac{b(\alpha(t) -0)-b(\alpha(t) +0)}{2i} \,h_p(y)   &
 a(\alpha(t)\! + \!0)\nu_p(y)\!+\!a(\alpha(t)\! - \!0)(1\!-\!\nu_p(y))\\
   \end{array}%
\!\!\!\right) %\label{eqn166}
\end{align*}
is invertible for every $(t,y)\in \tp\times \overline{\sR}$ and the
function
 \begin{align*}%\label{eqn17}
\smb (T(a)+H_{\alpha}(b))(t,y):=& \, a(t+0)\nu_p(y)+a(t-0)(1-\nu_p(y))\\
 & +
\mu(t) \, \frac{b(t+0)-b(t-0)}{2} \,h_p(y)
\end{align*}
where
 $$
\mu(t)=\left \{
 \begin{array}{rl}
1 &\text{if} \quad  t=t_+\\
-1 &\text{if}\quad t=t_-
 \end{array}
 \right . ,
 $$
does not vanish on $\{t_+,t_-\}\times \overline{\sR}$.

 \end{thm}
Note that although this result cannot be found in the literature, it
is not entirely new. It can be proved similarly to \cite{RS1990} by
using localization technique and the two-projection theorem. For
classical Toeplitz plus Hankel operators an analogous result is
presented in \cite{RSS:2011}. On the other hand, an index formula
can be established following ideas of \cite{RS:2012}. Moreover, if
the generating functions constitute a Fredholm matching pair
$(a,b)$, the kernel and cokernel of the operator
$T(a)+H_{\alpha}(b)$ can be described using results of
Section~\ref{s4}. For non-Fredholm matching pairs $(a,b)$ the
situation is more complicated. However, if $a$ and $b$ are piecewise
continuous functions, this case is still treatable. For example, the
following result is true.

 \begin{thm}\label{tt6}
Let $a,b\in PC$ and $(a,b)$ be a matching pair. If the operator
$T(a)+H_{\alpha}(b):H^p\to H^p$ is Fredholm, then there is an
interval $(p,p_0)$, $p<p_0$ such that for all $r\in(p,p_0)$ the pair
$(a,b)$ and both operators $T(a)\pm H_{\alpha}(b):H^r\to H^r$ are
Fredholm,
\begin{equation*}%\label{eq7.1}
 \begin{aligned}%\label{eq}
\ker (T(a)+H_{\alpha}(b))\left |_{H^r\to H^r} \right .&
 =\ker( T(a)+H_{\alpha}(b))\left |_{H^p\to H^p}\,, \right .\\
\coker (T(a)+H_{\alpha}(b))\left |_{H^r\to H^r} \right .&
 =\coker (T(a)+H_{\alpha}(b))\left |_{H^p\to H^p}\, , \right .\\
\end{aligned}
\end{equation*}
and the kernel and cokernel of the operator
$T(a)+H_{\alpha}(b):H^r\to H^r$ are described by Theorems
\ref{t5.1}--\ref{t4}, \ref{t5}.
\end{thm}
The proof is similar to the corresponding considerations of
\cite{DS:2014}.

 \section{How to determine the $\alpha$-factorization signature
for piecewise continuous functions}
 Suppose that a function $g\in PC$ satisfies the following two
conditions.
\begin{enumerate}
    \item $g g_\alpha=1$.
 \item The operator $T(g)$ is Fredholm on $H^p$.
 \end{enumerate}
In order to determine $\boldsymbol\sigma_{\!\alpha}( g)$ we need a
special factorization of the function $g$. More precisely, this
function has to be represented as
   \begin{equation}\label{eq21}
 g= \psi_{\beta_+,t_+} g_+, \, \text{ or }\,  g=\psi_{\beta_-,t_-} g_-,
 \end{equation}
where the function $\psi_{\beta_{\pm},t_{\pm}}$ and $g_{\pm}$
possess the following properties.
 \begin{enumerate}
\item The function $\psi_{\beta_+,t_+}$ has
a jump at the point $t_+$ and  is continuous on the arc
$\sT\setminus \{t_+\}$.

\item The function $\psi_{\beta_-,t_-}$ has
a jump at the point $t_-$ and  is continuous on the arc
$\sT\setminus\{t_-\}$.

\item  $\psi_{\beta_+,t_+} ( \psi_{\beta_+,t_+})_\alpha=1, \quad
 \psi_{\beta_-,t_-} ( \psi_{\beta_-,t_-})_\alpha=1.$

\item The function $g_+$ and $g_-$ are continuous at the
points $t_+$ and $t_-$, respectively.
 \end{enumerate}
It is clear that if such factorizations take place, then $g_+,g_-\in
PC$ and
 $$
g_+ (g_+)_\alpha =1, \quad g_- (g_-)_\alpha =1
 $$
It turns out that all the factorizations mentioned do really exist,
and below we show how to construct them.

Let  $z\neq 0$ be a complex number, and let $\arg z$ stand for that
value of the argument of $z$, which is located in the interval
$(-\pi, \pi]$. For $\beta\in\sC$, $\re \beta\in (-1/q,1/p)$, and
$\tau\in \sT$ consider the function $\vp_{\beta, \tau}(t)\in PC$
defined by
\begin{equation}\label{Eq50}
\vp_{\beta, \tau}(t):=\exp\{i\beta \arg (-t/\tau) \}, \quad t\in\sT.
\end{equation}
It is easily seen that $\vp_{\beta, \tau}$ has at most one
discontinuity, namely, a jump at the point $\tau$ and
\begin{equation*}
\vp_{\beta, \tau}(\tau+0)=\exp\{-i\pi\beta\}, \quad \vp_{\beta,
\tau}(\tau-0)=\exp\{i\pi\beta\}.
\end{equation*}
Recall a useful factorization of the function $\vp_{\beta,\tau}$,
viz.
\begin{equation}\label{eqN2}
 \vp_{\beta,\tau}(t)=\xi_{-\beta}(t)\eta_\beta(t),
\end{equation}
where
 \begin{align*}
\xi_\beta(t) =\left ( 1-\frac{\tau}{t} \right )^\beta & :=\exp\left
\{ \beta
  \log \left | 1-\frac{\tau}{t} \right | + i \beta \arg  \left ( 1-\frac{\tau}{t} \right )
  \right \}, \\
  \eta_\beta(t)=\left (1-\frac{t}{\tau}\right )^\beta & :=\exp\left \{ \beta
  \log  \left |1-\frac{t}{\tau} \right | + i \beta \arg \left ( 1-\frac{t}{\tau} \right )
  \right \}.
\end{align*}
Note that the representation \eqref{eqN2} is a Wiener-Hopf
factorization with the factorization index zero. Therefore, the
Toeplitz operator with generating function \eqref{Eq50} is
invertible on the space $H^p$ (see \cite[Sections 5.35 and
5.36]{BS:2006}).

Let $t_+$ and $T_-$ be the fixed  points of the mapping $\alpha$.
Recall that $t_+=(1+ \lambda)/\overline{\beta}$, $t_-=(1-
\lambda)/\overline{\beta}$, where $\lambda:=i \sqrt{|\beta|^2-1}$.
According to \cite[Sections 5.35 and 5.36]{BS:2006}, a function
$g\in PC$, such that $T(g)$ is Fredholm, can be written in one of
the following form
 \begin{equation}\label{rep}
 g= \vp_{\beta_+,t_+} g_1, \quad g= \vp_{\beta_-,t_-} g_2
 \end{equation}
where $\re \beta_{\pm,t_{\pm}}\in (-1/q,1/p)$ and the functions
$g_1,g_2\in PC$ are continuous at the points $t_+$ and $t_-$,
respectively. However,  the representations \eqref{rep} are of no
use in the present situation since the products
$\vp_{\beta_{\pm,t_{\pm}}}(\vp_{\beta_{\pm,t_{\pm}}})_{\alpha}$ are
not equal to $1$. Inspired  by  representation \eqref{eqN2} and by
Theorem \ref{t3} we define functions $\psi_{\beta_+,t_+}$ and
$\psi_{\beta_-,t_-}$ by
 \begin{align}\label{eq51}
 \psi_{\beta_+,t_+}(t)&:=
 \eta_{\beta_+,t_+}(t)\eta_{-\beta_+,t_+}(\alpha(t)), \quad t\neq
 t_+\,, \\
 \psi_{\beta_-,t_-}(t)&:=
 \eta_{\beta_-,t_-}(t)\eta_{-\beta_-,t_-}(\alpha(t)), \quad t\neq
 t_-\,. \label{eq52}
\end{align}
We are going to study properties of the functions \eqref{eq51} and
\eqref{eq52}. Let us deal with the function \eqref{eq51}. The other
one can be treated analogously. First of all, we note that
  $$
\eta_{\beta_+,t_+}\in H^q, \quad \eta_{\beta_+,t_+}^{-1}\in H^p,
\quad (\eta_{-\beta_+,t_+})_\alpha\in \overline{H^p}, \quad
(\eta_{-\beta_+,t_+})_\alpha^{-1}\in \overline{H^q},
  $$
where $1/p+1/q=1$. Thus if one shows that $\psi_{\beta_+,t_+}$ is a
piecewise continuous function, then representation \eqref{eq51} is a
weak Wiener--Hopf factorization for this function.

Due to our agreement about the choice of $\arg z$, $z\neq 0$, one
has
 $$
\arg\frac{t-t_+}{ \alpha(t)-t_+ } + 2k(t) \pi=\arg \left
(1-\frac{t}{t_+} \right ) -\arg \left (1-\frac{\alpha(t)}{t_+}
\right ) ,
 $$
where $k(t)$ is a uniquely determined integer, $t\neq t_+$. Now we
note that
 $$
\alpha(t)-t_+=\alpha(t)-\alpha(t_+)=\frac{(|\beta|^2
-1)(t-t_+)}{(\overline{\beta}t -1)(\overline{\beta}t_+ -1)}.
 $$
Therefore,
 \begin{equation}\label{Eq53}
 \frac{t-t_+}{\alpha(t)-t_+}=
  \frac{(\overline{\beta}t-1)i \sqrt{|\beta|^2-1}}{|\beta|^2-1},
  \quad t\neq t_+,
 \end{equation}
and
 \begin{equation*}%\label{Eq54}
\lim_{t\to t_+} \frac{t-t_+}{\alpha(t)-t_+}=-1.
 \end{equation*}
Consider now the argument of the function \eqref{Eq53}. It is equal
to $\arg(i(\overline{\beta}t-1))$, and we will show that
 $$
\lim_{t\to t_+\pm 0}\arg(i(\overline{\beta}t-1)) =\mp\pi.
 $$
In order to study $\arg(i(\overline{\beta}t-1))$ represent function
$t\mapsto i(\overline{\beta}t-1)$ as
 $$
i(\overline{\beta}t-1)= i(\overline{\beta}t_+\exp(-i\theta)-1),
\quad \theta \in (0,2\pi).
 $$
Using the relation $\exp(-i\theta)= \cos \theta - i\sin\theta$, one
can rewrite it in the form
  \begin{equation}\label{EqR}
 i(\overline{\beta}t-1)= \sqrt{|\beta|^2-1}\,(\sin\theta- \cos
\theta) + i(\cos \theta+\sin \theta - 1).
 \end{equation}
Consequently,
\begin{align*}%\label{Eq54}
&\lim_{t\to t_+-0} \arg \left ( \frac{t-t_+}{\alpha(t)-t_+}\right
)=\lim_{\theta\to0}\arg(i(\overline{\beta}t-1) )\\
 &=\lim_{\theta\to 0}\arg( \sqrt{|\beta|^2-1}\,(\sin\theta- \cos
\theta) + i(\cos \theta+\sin \theta - 1))=\pi ,
\end{align*}
because $\cos\theta +\sin\theta-1>0$ for $\theta$ small enough.

Analogously,
  \begin{equation*}
  \lim_{t\to t_++0} \arg \left ( \frac{t-t_+}{\alpha(t)-t_+}\right
)=-\pi.
 \end{equation*}
Indeed, using the relation \eqref{EqR} once again, we obtain
\begin{align*}%\label{Eq}
\lim_{t\to t_++0} \arg \left ( \frac{t-t_+}{\alpha(t)-t_+}\right )
&=\lim_{\theta\to 2\pi}\arg( \sqrt{|\beta|^2-1}\,(\sin\theta- \cos
\theta)\\ &\quad + i(\cos \theta+\sin \theta - 1))=-\pi.
\end{align*}
It is clear that the function
  $$
t \to \arg \left ( \frac{t-t_+}{\alpha(t)-t_+}\right )
  $$
is continuous on $\sT\setminus \{t_+\}$. Hence, it is piecewise
continuous on $\sT$, having only one jump at the point $t_+$ if
$\beta\neq 0$. Moreover, since both functions $\psi_{\beta_+,t_+}$
and
 $$
\exp\left \{ \beta_+ \log \left | \frac{t-t_+}{\alpha(t)-t_+} \right
| + i \beta_+ \arg  \left (\frac{t-t_+}{\alpha(t)-t_+} \right )
  \right \},
 $$
are continuous on  $\sT\setminus \{t_+\}$ and  do not vanish there,
we obtain that $t \to k(t)$, $t\neq t_+$ is also continuous on
$\sT\setminus \{t_+\}$, hence it is even constant. Now let us choose
$t_0\in\sT$ such that $\alpha(t_0)=t_+$. Then
 $$
\arg\left ( 1-\frac{\alpha(t_0)}{t_+}\right )=0, \quad \arg\left (
1-\frac{t_0}{t_+}\right )\in (-\pi,\pi],
 $$
so $k(t)=0$.

 Thus it is shown that
  $$
\arg \left (1-\frac{t}{t_+} \right ) -\arg \left
(1-\frac{\alpha(t)}{t_+} \right )= \arg\frac{t-t_+}{ \alpha(t)-t_+ }
, \quad t\neq t_+.
 $$
As a consequence,  the function $\psi_{\beta_+,t_+}$ possesses the
following properties.

 \begin{enumerate}

 \item The function $\psi_{\beta_+,t_+}$ has a jump at the point
 $t_+$,
 \begin{align*}%\label{Eq}
    \lim_{t\to t_+ -0}  \psi_{\beta_+,t_+}(t)&= \exp(i \pi
\beta_+),\\
   \lim_{t\to t_++0} \psi_{\beta_+,t_+}(t)&= \exp(- i \pi
\beta_+),
\end{align*}
and it is continuous on $\sT\setminus \{t_+\}$.

 \item The function $\psi_{\beta_+, t_+}$ satisfies the relation
 $$
\psi_{\beta_+, t_+} (\psi_{\beta_+, t_+})_\alpha=1.
 $$
 \end{enumerate}

Now let us show  that the Toeplitz operator $T(\psi_{\beta_+,t_+})$
is invertible in $H^p$. Indeed, one has
 $$
\psi_{\beta_+,t_+}=\vp_{\beta_+,t_+} h, \quad h=
\frac{\psi_{\beta_+,t_+}}{\vp_{\beta_+,t_+}},
 $$
where $\vp_{\beta_+,t_+}$ is defined by \eqref{Eq50} and $h$ is a
continuous function which does not vanish on $\sT$. Then the
operator $T(\psi_{\beta_+,t_+})- T(\vp_{\beta_+,t_+}) T(h)$ is
compact. Since $T(\vp_{\beta_+,t_+}) T(h)$ is a Fredholm operator,
the operator $T(\psi_{\beta_+,t_+})$ is also Fredholm, and the
factorization \eqref{eq51} is in fact a Wiener--Hopf factorization
with the factorization index $0$. Hence, the operator
$T(\psi_{\beta_+,t_+})$ is invertible in $H^p$.

Note that the factorization \eqref{eq51} is not normalized, that is
$\eta_{-\beta_+, t_+}$ is not equal to $1$ at infinity. However, the
normalization can be easily achieved by multiplying the factor
$(1-\alpha(t)/t_+)^{\beta_+}$ by the number
$c:=(1-(\overline{\beta})^{-1}/t_+)^{-\beta_+}$ and the factor
$(1-\alpha(t)/t_+)^{-\beta_+}$ by the number $c^{-1}$.

Thus the factorization signature
$\boldsymbol\sigma(\psi_{\beta_+,t_+})$ of the function
$\psi_{\beta_+,t_+}$ is equal to $1$. This can be also obtained by
observing that
 $$
\psi_{\beta_+,t_+}(t_-)=1.
 $$

   \begin{thm}\label{t7}
Let $g_1\in PC$ be as above. Then
 $$
\boldsymbol\sigma_{\!\alpha}( g)=\boldsymbol\sigma_{\!\alpha}(
g_+)=\boldsymbol\sigma_{\!\alpha}( g_-).
 $$
  \end{thm}
The proof of these results run in parallel to the corresponding
results of \cite[Section 8]{DS:2014}.

Thus if the generating functions $a$ and $b$ of the generalized
Toeplitz plus Hankel operator $T(a)+H_{\alpha}(b)$ are piecewise
continuous and satisfy the matching condition \eqref{eq15}, then
representations \eqref{eq21} allows one to find $\alpha$-signature
of the corresponding auxiliary functions $c, \overline{c}$ and $d,
\overline{d}$ and, consequently, to obtain an effective and complete
description of the kernel and cokernel of the operator under
consideration.

  \begin{rem}
If $a\in L^\infty$ is a matching function having one-sided limits at
the points $t=t_{\pm}$, then Theorem \ref{t7} is still true.
Moreover, the $\alpha$-factorization signature can be effectively
calculated even in the case where $a$ has one-sided limits only at
the one of the fixed points $t_+$ or $t_-$.
  \end{rem}

%\newpage
%  \bibliographystyle{acm}
 %\bibliography{E:/TeXX/JabRef/DWP_2010}

\begin{thebibliography}{10}

\bibitem{BE2004}
{\sc Basor, E.~L., and Ehrhardt, T.}
\newblock Factorization theory for a class of {T}oeplitz {$+$} {H}ankel
  operators.
\newblock {\em J. Operator Theory 51}, 2 (2004), 411--433.

\bibitem{BE2006}
{\sc Basor, E.~L., and Ehrhardt, T.}
\newblock Factorization of a class of {T}oeplitz + {H}ankel operators and the
  {$A_p$}-condition.
\newblock {\em J. Operator Theory 55}, 2 (2006), 269--283.

\bibitem{BE:2013}
{\sc Basor, E.~L., and Ehrhardt, T.}
\newblock Fredholm and invertibility theory for a special class of {T}oeplitz +
  {H}ankel operators.
\newblock {\em J. Spectral Theory 3}, 3 (2013), 171--214.

\bibitem{BS:2006}
{\sc B{\"o}ttcher, A., and Silbermann, B.}
\newblock {\em Analysis of {T}oeplitz operators}, second~ed.
\newblock Springer Monographs in Mathematics. Springer-Verlag, Berlin, 2006.
%\newblock Prepared jointly with Alexei Karlovich.

\bibitem{DS:2013}
{\sc Didenko, V.~D., and Silbermann, B.}
\newblock Index calculation for {T}oeplitz plus {H}ankel operators with
  piecewise quasi-continuous generating functions.
\newblock {\em Bull. London Math. Soc. 45}, 3 (2013), 633--650.

\bibitem{DS:2014b}
{\sc Didenko, V.~D., and Silbermann, B.}
\newblock The {C}oburn-{S}imonenko {T}heorem for some classes of
  {W}iener--{H}opf plus {H}ankel operators.
\newblock {\em Publications de l'Institut Mathe'matique 96(100)} (2014), 85--102.

\bibitem{DS:2014a}
{\sc Didenko, V.~D., and Silbermann, B.}
\newblock Some results on the invertibility of {T}oeplitz plus {H}ankel
  operators.
\newblock {\em Annalas Academie Scientarium Fennicae, Mathematica 39\/} (2014),
  439--446.

\bibitem{DS:2014}
{\sc Didenko, V.~D., and Silbermann, B.}
\newblock Structure of kernels and cokernels of {T}oeplitz plus {H}ankel
  operators.
\newblock {\em Integral Equations Operator Theory 80}, 1 (2014), 1--31.

\bibitem{Ehr:2004h}
{\sc Ehrhardt, T.}
\newblock {\em Factorization theory for {T}oeplitz+{H}ankel operators and
  singular integral operators with flip}.
\newblock Habilitation {T}hesis, Technische Universit{\"a}t Chemnitz, 2004.

\bibitem{KLR:2007}
{\sc Kravchenko, V.~G., Lebre, A.~B., and Rodr{\'{\i}}guez, J.~S.}
\newblock Factorization of singular integral operators with a {C}arleman shift
  via factorization of matrix functions: the anticommutative case.
\newblock {\em Math. Nachr. 280}, 9-10 (2007), 1157--1175.

\bibitem{KLR:2009}
{\sc Kravchenko, V.~G., Lebre, A.~B., and Rodr{\'{\i}}guez, J.~S.}
\newblock Factorization of singular integral operators with a {C}arleman
  backward shift: the case of bounded measurable coefficients.
\newblock {\em J. Anal. Math. 107\/} (2009), 1--37.

\bibitem{LS1987}
{\sc Litvinchuk, G.~S., and Spitkovskii, I.~M.}
\newblock {\em Factorization of measurable matrix functions}, vol.~25 of {\em
  Operator Theory: Advances and Applications}.
\newblock Birkh\"auser Verlag, Basel, 1987.
%\newblock Translated from the Russian by Bernd Luderer, With a foreword by
  %Bernd Silbermann.

\bibitem{RSS:2011}
{\sc Roch, S., Santos, P.~A., and Silbermann, B.}
\newblock {\em Non-commutative {G}elfand theories}.
\newblock Universitext. Springer-Verlag London Ltd., London, 2011.
\newblock {A tool-kit for operator theorists and numerical analysts}.

\bibitem{RS1990}
{\sc Roch, S., and Silbermann, B.}
\newblock {\em Algebras of convolution operators and their image in the
  {C}alkin algebra}, vol.~90 of {\em Report MATH}.
\newblock Akademie der Wissenschaften der DDR Karl-Weierstrass-Institut f\"ur
  Mathematik, Berlin, 1990.
%\newblock With a German summary.

\bibitem{RS:2012}
{\sc Roch, S., and Silbermann, B.}
\newblock {A handy formula for the Fredholm index of Toeplitz plus Hankel
  operators}.
\newblock {\em Indag. Math. 23}, 4 (2012), 663--689.

\bibitem{Si:1987}
{\sc Silbermann, B.}
\newblock The {$C^*$}-algebra generated by {T}oeplitz and {H}ankel operators
  with piecewise quasicontinuous symbols.
\newblock {\em Integral Equations Operator Theory 10}, 5 (1987), 730--738.
\newblock Toeplitz lectures 1987 (Tel-Aviv, 1987).

\end{thebibliography}

 \end{document}